\documentclass[12pt]{article}

\usepackage{import}
\usepackage[left=3cm,right=3cm,top=2.8cm,bottom=3cm, marginparsep=3mm]{geometry}
\usepackage{changepage}
\usepackage{adjustbox}

\usepackage{graphicx} 
\usepackage[dvipsnames]{xcolor}

\usepackage{amssymb}
\usepackage{amsfonts}
\usepackage{amsmath}
\usepackage{amscd}
\usepackage{aligned-overset}
\usepackage{easy-todo}
\usepackage{enumitem}
\usepackage{caption}
\usepackage{changepage}

\usepackage{amsthm}

\theoremstyle{plain}
\newtheorem{theorem}{Theorem}
\newtheorem{proposition}[theorem]{Proposition}

\newtheorem{conjecture}{Conjecture}
\newtheorem*{conjecture*}{Conjecture}
\newtheorem{lemma}[theorem]{Lemma}
\newtheorem{corollary}[theorem]{Corollary}

\theoremstyle{definition}
\newtheorem{definition}[theorem]{Definition}
\newtheorem{example}[theorem]{Example}
\newtheorem{remark}[theorem]{Remark}

\usepackage{tikz-cd}
\usepackage{tikz}
\usetikzlibrary{shapes, shapes.geometric, arrows, arrows.meta, decorations.markings, decorations.pathreplacing, decorations.pathmorphing}
\pgfdeclarelayer{background} 
\pgfdeclarelayer{foreground}
\pgfsetlayers{background,main,foreground}

\newcommand{\N}{\mathbb{N}}
\newcommand{\Z}{\mathbb{Z}}

\newcommand{\R}{\mathbb{R}}

\newcommand{\Sb}{\mathbb S}

\newcommand{\B}{\mathrm{B}}

\renewcommand{\Im}{\mathrm{Im}}

\newcommand{\Stab}{\mathrm{Stab}}

\newcommand{\G}{\mathcal{G}} 
\newcommand{\Gbar}{\mathcal{G}_\mathbb{\R}} 
\newcommand{\Hc}{\mathcal{H}} 
\newcommand{\Cc}{\mathcal C} 
\newcommand{\Ec}{\mathcal E} 
\newcommand{\Ac}{\mathcal A}

\DeclareMathOperator{\Conv}{CvxHul}
\newcommand{\E}{\mathrm{E}}
\newcommand{\F}{\mathrm{F}}
\DeclareMathOperator{\Cl}{Cl} 
\DeclareMathOperator{\Iso}{Iso}

\newcommand{\say}[1]{``#1"}
\newcommand{\Stoll}{\mathrm{Stoll}}
\newcommand{\Geo}{\mathrm{Geo}}

\renewcommand{\le}{\leqslant}
\renewcommand{\ge}{\geqslant}

\newcommand{\onto}{\twoheadrightarrow}
\newcommand{\longto}{\longrightarrow}
\newcommand{\acts}{\curvearrowright} 

\newcommand{\la}{\left\langle}
\newcommand{\ra}{\right\rangle}

\newcommand{\abs}[1]{\left| #1 \right|}
\newcommand{\norm}[1]{\left\|#1\right\|}

\DeclareMathOperator{\id}{id}

\DeclareMathOperator{\newdiff}{d} 
\newcommand{\dif}{\newdiff\!}

\usepackage[backend=bibtex, style=numeric, sorting=nyt, 
    doi=false,isbn=false,url=false,eprint=false, maxnames=50]{biblatex}
\bibliography{horofunction}

\usepackage{hyperref}

\title{Horofunctions on the Heisenberg and Cartan groups}
\author{Corentin Bodart and Kenshiro Tashiro}
\date{\today}
\begin{document}

\maketitle
\vspace*{-2mm}

\begin{abstract}
    We study the horofunction boundary of finitely generated nilpotent groups, and the natural group action on it. More specifically, we prove the followings:
    \begin{itemize}[leftmargin=6mm, rightmargin=2mm]
        \item For discrete Heisenberg groups, we classify the orbits of Busemann points. As a byproduct, we observe that the set of orbits is finite and the set of Busemann points is countable.
        Furthermore,
        by using the approximation with Lie groups, we also observe that the entire horoboundary is uncountable.
        
        \item For the discrete Cartan group, we exhibit an continuum of Busemann points, disproving a conjecture of Tointon and Yadin. As a byproduct, we prove that the group acts non-trivially on its reduced horoboundary, disproving a conjecture of Bader and Finkelshtein.
    \end{itemize}
\end{abstract}

\bigbreak
\renewcommand*{\thetheorem}{\Alph{theorem}}
The horofunction boundary (or horoboundary) was introduced by Gromov in \cite{Gromov}, as a way to compactify proper metric spaces, e.g., Cayley graphs. His construction was inspired by the Busemann function compactification of Hadamard spaces (i.e., complete CAT(0) spaces). Both compactifications share nice properties, for instance isometric actions on the metric space can be extended to actions on the boundary. This fact induces a natural group action of any finitely generated group on its horoboundary via the isometric action on its Cayley graph. Moreover, every (almost) geodesic ray converges to a unique horofunction. Such limits are called Busemann point. However, outside of the CAT(0) setting, horofunctions which are not Busemann points may (and often do \cite{webster-winchester2006}) exist. Therefore, we get two actions of geometric origins to satisfy our curiosity (and try to grasp insights into our groups): the action on the whole horoboundary, and on the subset of Busemann points.

\bigbreak

In this article, our focus lies on horofunction boundaries of finitely generated nilpotent groups. The horoboundary of nilpotent groups has received its fair share of attention, starting with abelian groups \cite{Rieffel,develin2002,webster-winchester2006}. One of the primary motivations is the following observation,
due to Anders Karlsson \cite{lecture_notes}:
\smallbreak
\begin{adjustwidth}{4mm}{4mm}
If the horoboundary contains a finite orbit $\G\cdot \varphi$, then $\G$ contains a finite-index subgroup which factors onto $\Z$. ($\varphi\colon \Stab(\varphi)\to \Z$ is a non-trivial homomorphism.)
\end{adjustwidth}
\medbreak
\noindent In particular, if one could prove that all finitely generated groups of polynomial growth admit a finite orbit in their horoboundary,
this would provide an alternative proof of Gromov’s theorem on groups of polynomial growth \cite{Gromov_poly}.

\newpage

The existence of finite orbits was indeed proved for nilpotent groups \cite{walsh2011}. In the original spirit, this was checked for groups of linear growth by Tointon and Yadin \cite{TointonYadin} using only growth considerations. One way to ensure the existence of finite orbits is to check that the boundary itself is small. Formally, Karlsson proved that, if $\partial(\G,S)$ is countable, then $\G\acts\partial(\G,S)$ admits a finite orbit \cite[Corollary 5]{lecture_notes}. In this direction, Tointon and Yadin conjectured that the horoboundary of groups of polynomial growth should contain only countably many Busemann points.

\bigbreak

More specifically, we consider the horofunctions and Busemann points of (lattices of) the $(2k+1)$-dimensional Heisenberg groups $\Hc_k$ and the Cartan group $\Cc$. These are prototypical examples of $2$-step and $3$-step nilpotent groups respectively. 

\bigbreak

For the Heisenberg groups and its generalization, we shall give a complete classification of the orbits of Busemann points (see Theorem \ref{theorem:nheisenberg orbits} for a precise statement).
\begin{theorem}
    Let $\Hc$ be a finitely generated, torsionfree $2$-step nilpotent group with $[\Hc,\Hc]\simeq \Z$, and $S$ any generating set. Then $\partial(\Hc,S)$ contains only finitely many orbits of Busemann points. In particular, there are only countably many Busemann points.
\end{theorem}
\noindent On the other hand, the entire horoboundary can already be large:
\begin{theorem}[Corollary \ref{cor:Heisenberg}]
    Let $\Hc_1$ be the $3$-dimensional Heisenberg group, and $S$ any generating set. Then $\partial(\Hc_1,S)$ has the cardinality of the continuum.
\end{theorem}
\noindent This indicates that \say{almost all} horofunctions of (non virtually abelian) nilpotent groups may not be Busemann points.

\bigbreak

For the Cartan group $\Cc=\la x,y \mid \text{3-step nilpotent}\ra$ (with generating set $S=\{x,y\}^\pm$), we use a geometrical model to construct an explicit family of geodesic rays $(\gamma_u)_u$, parametrized by unit vectors $u\in\Sb^1$, with the following property:
\begin{theorem}[Theorem \ref{thm:cartan}]The Busemann points $b_{\gamma_u}$ are pairwise distinct, and their orbits in  the reduced horoboundary $\partial^r(\Cc,S)$ are infinite.
\end{theorem}
\noindent Here $\partial^r (\Cc,S)\coloneqq \partial(\Cc,S)/C_b(\Cc)$ ($C_b(\Cc)$ is the space of bounded functions on $\Cc$). This disproves the conjecture of Tointon and Yadin, as well as a conjecture of Bader and Finkelshtein \cite{Bader} (\say{$\G\acts \partial^r(\G,S)$ is trivial for all finitely generated nilpotent groups $\G$}).

\bigbreak

Based on our observations, it is reasonable to refine the first conjecture as follows:
\begin{conjecture} Let $\G$ be a group and $S$ a finite generating set.
    \begin{enumerate}[leftmargin=8mm, label={\normalfont(\alph*)}]
        \item The horoboundary $\partial(\G,S)$ is countable if and only if $\G$ is virtually abelian.
        \item If $\G$ is virtually $2$-step nilpotent, then the set of orbits of Busemann points in $\partial(\G,S)$ is finite. In particular, there are at most countably many Busemann points.
    \end{enumerate}
\end{conjecture}
\noindent We do not adventure in a complete characterization for part (b), as many other groups (or rather marked groups) may enter the picture. For instance, filiform groups (such as the Engel group), or pairs $(\G,S)$ with intermediate geodesic growth (including virtually Cartan groups! \cite{Intermediate}) may only have countably many Busemann points.

\bigbreak

\noindent \textbf{Acknowledgements.} The first named author was supposed by the Swiss NSF grant 200020-200400. The second named author is partially supported by JSPS KAKENHI grant number 23K03104. 

\renewcommand*{\thetheorem}{\arabic{theorem}}
\counterwithin{theorem}{section}
\section{Preliminaries}

\subsection{Horofunctions and Busemann points}

We first define the horofunction boundary of a general proper metric space $(X,d)$. In what follows, we will choose $(X,d)$ to be a finitely generated group with a word metric, or a simply connected nilpotent Lie group with a sub-Finsler metric.
\begin{definition}[Horofunction boundary] Let $(X,d)$ be a proper metric space, together with a base point $e\in X$. Consider the embedding $\iota\colon X\to C(X,\R)$ sending
\[ x \mapsto \big(\varphi_x(\ast)=d(x,\ast)-d(x,e)\big), \]
where $C(X,\R)$ is the set of continuous functions $X\to\R$ with the topology of uniform convergence on compact sets. The \emph{horofunction boundary} (or \emph{horoboundary}) is defined as $\partial(X,d)\coloneqq \Cl(\iota(X))\setminus \iota(X)$, and functions $\varphi\in\partial(X,d)$ are called \emph{horofunctions}.
\end{definition}

\bigbreak

We observe that, since $X$ is proper (hence second countable), every horofunction can be seen as a limit $\varphi(\ast)=\lim_{n\to\infty}\varphi_{x_n}(\ast)$ for some sequence $(x_n)\subseteq X$.

\begin{definition}
Let $T\subseteq [0,\infty)$ be an unbounded subset (e.g. $T=\N$),
and $\gamma:T\to X$ be a mapping. We say that
\begin{itemize}[leftmargin=8mm]
    \item $\gamma$ is a \emph{geodesic ray} if, for all $s,t\in T$, we have $d(\gamma(s),\gamma(t))=|t-s|$.
    \item $\gamma$ is an \emph{almost geodesic ray} if for all $\varepsilon>0$,
    there is $N\in\N$ such that
    \[ \forall s,t\in T\cap [N,\infty),\quad |d(\gamma(t),\gamma(s))+d(\gamma(s),\gamma(0))-t|<\varepsilon.\]
\end{itemize}
A horofunction is called a \emph{Busemann point} if it is the limit of an almost geodesic ray $\gamma:T\to X$. We denote the associated Busemann point by $b_\gamma\colon X\to\R$.
\end{definition}
Let $g$ be an isometry of $(X,d)$.
By \cite[Proposition 4.10]{Rieffel}, $g$ extends uniquely to an homeomorphism of $\mathrm{Cl}(\iota(X))$, sending $\partial(X,d)$ to itself. This defines a (left) group action  $\Iso(X)\curvearrowright \partial (X,d)$ by homeomorphisms, given by the explicit formula
\[ (g\cdot \phi)(x)= \phi\big(g^{-1}(x)\big) -\phi\big(g^{-1}(e)\big).\]
Moreover, the action $\Iso(X)\acts (X,d)$ sends almost geodesic rays to almost geodesic rays, hence $\Iso(X)\acts\partial(X,d)$ sends Busemann points to Busemann points.

\subsection{Reduced horoboundary and rough isometries} \label{ssec:reduced}

In general it is difficult to compute horofunctions since this requires to compute some distances \emph{exactly}. In particular, the horoboundary of a Cayley graph (i.e., its topology) typically depends on the choice of a generating set. However,
once we ignore bounded errors, working with horofunctions often become much more manageable. In some cases, the resulting boundary even becomes a quasi-isometry invariant.

\medbreak

\begin{definition}
    Let $C_b(X,\R)$ be the set of bounded continuous functions $X\to\R$.
    The quotient space $\partial^r(X,d)=\partial(X,d)/C_b(X,\R)$ is called the \emph{reduced horoboundary} of $(X,d)$. Elements of the reduced horoboundary are denoted by $[\varphi]$.
\end{definition}
\begin{example}
    The reduced horoboundary of a proper Gromov hyperbolic space is homeomorphic to its Gromov boundary \cite{Webster-winchester-hyp-2005}. It is a quasi-isometry invariant.
\end{example}

\bigbreak

Since the reduced horoboundary ignores additive errors, it fits well with the coarse setting. The coarse analog of isometries is the notion of rough isometry. We will see that reduced horoboundaries are stable under the rough isometry:
\begin{definition}
A map $f\colon (X,d_X)\to (Y,d_Y)$ is a \emph{rough isometry} if
    \begin{enumerate}[leftmargin=8mm, label=(\alph*)]
        \item there exists $C_1\ge 0$ such that $\abs{d_Y\big(f(x),f(x')\big)-d_X(x,x')}\le C_1$ for all $x,x'\in X$,
        \item there exists $C_2\ge 0$ such that $d_Y\big(y,f(X)\big)\le C_2$ for all $y\in Y$.
    \end{enumerate}
    We also say $f$ is 
    a \emph{roughly isometric embedding} if it only satisfies condition (a).
\end{definition}

\begin{proposition}\label{proposition:equivalencehorofunction}
Let $(X,d_X),(Y,d_Y)$ be two proper metric spaces. Given a roughly isometric embedding $f\colon X\to Y$, we can define a surjective map
    \[ f^*\colon \begin{pmatrix} \Cl(\iota(f(X))) \cap \partial^r(Y,d_Y) & \longto & \partial^r(X,d_X) \\
    [\varphi] & \longmapsto & [\varphi\circ f] \end{pmatrix}. \]
    If $f$ is a rough isometry, then $f^*$ is a bijection $\partial^r(Y,d_Y)\to\partial^r(X,d_X)$.
\end{proposition}
\begin{proof} Fix $e\in X$ and $f(e)\in Y$ two basepoints.

\bigbreak
  
     \noindent $\blacktriangleright$ \emph{Well-defined:} Consider $\varphi\in \Cl(f(X))\cap\partial(Y,d_Y)$. We fix $(x_n)\subseteq X$ s.t.\ $f(x_n)\to \varphi$ in $\partial(Y,d_Y)$, and a subsequence $(x_{n_k})$ such that $x_{n_k}\to\psi\in \partial(X,d_X)$. We have
    \begin{align*}
        & \abs{\varphi(f(g))-\psi(g)} \\
        & \qquad = \lim_{k\to\infty} \abs{ \Big(d_Y\!\big(f(g),f(x_{n_k})\big)-d_Y\!\big(f(e),f(x_{n_k})\big)\Big) - \Big(d_X(g,x_{n_k})-d_X(e,x_{n_k})\Big) }\\
        & \qquad \le 2C_1
    \end{align*}
    so that $[\varphi\circ f]=[\psi]\in\partial^r(X,d_X)$. Moreover, given two representatives $\varphi,\varphi'\in[\varphi]$,
    \[ \norm{\varphi\circ f-\varphi'\circ f}_{L^\infty(X)} \le \norm{\varphi-\varphi'}_{L^\infty(Y)}<\infty \]
    so that $[\varphi\circ f]=[\varphi'\circ f]$.

    \bigbreak

    \noindent $\blacktriangleright$ \emph{Surjective:} Consider $[\psi]\in\partial^r(X,d_X)$, and $(x_n)\subseteq X$ s.t.\ $x_n\to\psi$. As $\Cl(f(X))\cap\partial(Y,d_Y)$ is compact, there exists $\varphi\in \Cl(f(X))\cap\partial(Y,d_Y)$ and a subsequence $(x_{n_k})$ such that $f(x_{n_k}) \to \varphi$. Using the first computation, we have $f^*([\varphi])=[\varphi\circ f]=[\psi]$.

    \bigbreak

    \noindent $\blacktriangleright$ \emph{New domain:} We suppose that $f$ is a rough isometry and prove that
    \[ \Cl(f(X)) \cap \partial^r(Y,d_Y)=\partial^r(Y,d_Y). \] Fix $[\varphi]\in \partial^r(Y,d_Y)$ and a sequence $(y_n)\subseteq Y$ such that $y_n\to \varphi$. As $f$ is a rough isometry, we can find $(x_n)\subseteq X$ such that $d_Y(y_n,f(x_n))\le C_2$. As $\partial(Y,d_Y)$ is compact, there is a subsequence $x_{n_k}$ such that $f(x_{n_k})\to \varphi'\in\partial(Y,d_Y)$. We have
    \begin{align*}
        &\abs{\varphi(g)-\varphi'(g)}\\
        & \qquad = \lim_{k\to\infty} \abs{ \Big(d_X(g,y_{n_k})-d_X\!\big(f(e),y_{n_k}\big)\Big) - \Big(d_Y\!\big(g,f(x_{n_k})\big)-d_Y\!\big(f(e),f(x_{n_k})\big)\Big) }\\
        & \qquad \le 2C_2
    \end{align*}
    so that $[\varphi]=[\varphi']\in \Cl(f(X))$.
    
    \bigbreak
    
    \noindent$\blacktriangleright$ \emph{Inverse:} There exists a rough isometry $g\colon Y\to X$ such that
    \[ d_X\!\big(g(f(x)) ,\id_X(x)\big)\le C_2. \]
    It follows that $\abs{\psi\circ g\circ f - \psi} \le C_2$
    for all $\psi\in\partial(X,d_X)$, as horofunctions are $1$-Lipschitz, hence $f^*(g^*([\psi])=[\psi\circ g\circ f]=[\psi]$, i.e., $f^*\circ g^*=\id_{\partial^r(X,d_X)}$. Similarly, $g^*\circ f^*=\id_{\partial^r(Y,d_Y)}$. We conclude that $f^*$ is a bijection. \qedhere

\end{proof}

\subsection{Horofunctions on Cayley graphs}

Let us consider the special case of horofunctions on Cayley graphs. The graph structure allows to simplify some definitions and problems. For instance,
\begin{itemize}[leftmargin=8mm]
    \item The topology is discrete, hence the topology of uniform convergence on compact sets coincides with the topology of pointwise convergence.
    \item Every Busemann point is the limit of a geodesic \cite[Lemma 2.1]{webster-winchester2006}.
\end{itemize}
We also have the following folklore lemma:
\begin{lemma}[{see eg.\ \cite[Lemma 2.4]{TointonYadin}}] \label{lem:from_e}
In a graph $(\G,d)$, any Busemann point $b_\gamma$ is represented by a geodesic $\eta$ starting at $\eta_0=e$.
\end{lemma}
\begin{proof}[Sketch of proof]
    The integer-valued sequence $d\big(\gamma_n,e\big)-d\big(\gamma_n,\gamma_0\big)$ is decreasing, so it is constant for $n\ge N$, for some $N$. In particular, we have
    \[ \forall n\ge N,\quad d(\gamma_n,e)=d(\gamma_n,\gamma_N)+d(\gamma_N,e). \]
    This means that $\eta=([\text{geodesic path from $e$ to $\gamma_N$}], \gamma_N,\gamma_{N+1},\ldots)$ is a geodesic ray, starting at $e$, which eventually coincides with $\gamma$. In particular, $b_\eta=b_\gamma$.
\end{proof}

For the study of Cayley graphs and their horoboundary, we have a rich structure given by graph labellings, left translations and homomorphisms. We consider $\G$ a group, with $S$ a finite symmetric generating set, and $d_S$ the induced word metric. Any infinite path $(x_n)_{n\ge 0}$ starting at $x_0=e$ can be parametrized by an infinite word $\gamma=s_1s_2\cdots \in S^\infty$, with $s_i=\gamma_i^{-1}\gamma_{i+1}$. We denote prefixes of $\gamma$ by $\gamma_n\coloneqq s_1s_2\cdots s_n\in S^n$. The associated group elements are denoted as $\bar{\gamma}_n\in \G$ (so $\bar\gamma_n=x_n$).

\begin{definition}
    The language of infinite geodesics is
    \[ \Geo(\G,S)=\{\gamma\in S^\infty\mid  \forall n,\; d_S(e,\bar\gamma_n)=n\}. \]
\end{definition}
In this formalism, Lemma \ref{lem:from_e} easily gives the following:
\begin{lemma} \label{lemma:action}
    In a Cayley graph $(\G,d_S)$,
    \begin{enumerate}[leftmargin=8mm, label={\normalfont(\alph*)}]
        \item Every Busemann point is represented by an infinite geodesic word $\gamma\in\Geo(\G,S)$.
        \item Two Busemann points lie in the same orbit if and only if they are represented by infinite geodesic words which are \emph{cofinal}, i.e., equal up to removing finite prefixes.
    \end{enumerate}
\end{lemma}
\begin{proof}
    Given $g\in\G$ and $\gamma\in\Geo(\G,S)$, the Busemann point $g\cdot b_\gamma$ is represented by a geodesic with the same labeling, starting from $g$ (instead of $e$). So the construction in Lemma \ref{lem:from_e} gives a geodesic word $\eta$ representing $g\cdot b_\gamma$ in the same cofinality class as $\gamma$.
    
    \medbreak
    
    \noindent Reciprocally, given two infinite geodesics $\gamma=uw$ and $\eta=vw$ (with $u,v\in S^*$ finite prefixes and $w\in S^\infty$ the common infinite suffix) we have $b_\gamma= \bar u \cdot b_w = \bar u\bar v^{-1}\cdot b_\eta$.
\end{proof}

The following proposition is helpful to distinguish/identify Busemann points in $\partial(\G,d)$ and $\partial^r(\G,d)$. The part (1) is essentially proven in \cite[Proposition 2.1]{walsh2011}
\begin{proposition} \label{prop:switch}
Consider two geodesic rays $\gamma$ and $\eta$ in a graph $(\G,d)$ starting at $\gamma_0=\eta_0=e$. The following assertions are equivalent:
\begin{itemize}[leftmargin=12mm]
    \item[{\normalfont(1a)}] $\gamma$ and $\eta$ converge to the same Busemann point, that is, $b_\gamma=b_\eta$ in $\partial(\G,d)$.
    \item[{\normalfont(1b)}] For all $n$, there exists $m\ge n$ such that $d(\gamma_m,\eta_n) = d(\eta_m,\gamma_n) = m-n$.
    \item[{\normalfont(1c)}] There is a geodesic having infinitely many points in common with both $\gamma$ and $\eta$.
\end{itemize}
In the coarse setting, we have the following equivalence:
\begin{itemize}[leftmargin=12mm]
    \item[{\normalfont(2a)}] $[b_\gamma]=[b_\eta]$ in $\partial^r(\G,d)$.
    \item[{\normalfont(2b)}] There exists a constant $C$ such that, for all $n$, there exists $m\ge n$ such that
    \[ d(\gamma_m,\eta_n),\; d(\eta_m,\gamma_n) \le  m-n + C. \]
\end{itemize}
\end{proposition}
\begin{proof}
    (2a) $\Rightarrow$ (2b): There exists $C\ge 0$ such that $\abs{b_\gamma(x)-b_\eta(x)}\le C$ for all $x\in\G$. In particular, this holds for $x=\gamma_n$. For $m\ge n$ large enough, we have
    \[ -C\le b_\gamma(\gamma_n)-b_\eta(\gamma_n) = \big(d(\gamma_m,\gamma_n) - d(\gamma_m,e) \big) - \big(d(\eta_m,\gamma_n) - d(\eta_m,e)\big) = m-n-d(\eta_m,\gamma_n) \]
    hence $d(\eta_m,\gamma_n)\le m-n+C$. Similarly, $d(\gamma_m,\eta_n)\le m-n+C$ for $m$ large enough.

    \medbreak

    \noindent (2b) $\Rightarrow$ (2a): Let $x\in\G$. For $m$ and $n$ well-chosen, we have
    \[ b_\eta(x)-b_\gamma(x) = \big(d(\eta_m,x)-m\big) - \big(d(\gamma_n,x)-n\big)\le d(\eta_m,\gamma_n)-(m-n) \le C. \]
    Similarly, we have $b_\gamma(x)-b_\eta(x)\le C$ for all $x\in\G$.
\end{proof}

\bigbreak

Let us consider a surjective homomorphism $f\colon\G\onto \Hc$. Then $T\coloneqq f(S)$ is a generating set of $\Hc$. We shall denote by $d_S,d_T$ the word metrics on $\G,\Hc$ respectively.

\begin{corollary} \label{lemma:quotient}
For each $t\in T$, fix an element $\tilde t\in S$ such that $\pi(\tilde t)=t$. For a set of infinite geodesic word $\Gamma\subseteq \Geo(\G,S)$, let $\partial(\Gamma) \subseteq \partial(\G,d_S)$ be the set of Busemann points induced from $\Gamma$, and $\partial^r(\Gamma)\subset \partial^r(\G,d_S)$ its quotient onto reduced horoboundary.
   \begin{enumerate}[leftmargin=8mm, label={\normalfont(\arabic*)}]
			\setcounter{enumi}{-1}
			\item The map
   \[ \pi^*\colon \begin{pmatrix} \Geo(\Hc,T)& \longto &  \Geo(\G,S)  \\ \gamma=t_1t_2\cdots & \longmapsto & \tilde\gamma=\tilde{t}_1\tilde{t}_2\cdots \end{pmatrix} \]
   is well-defined and injective.
\item The map
\[ \pi_*\colon \begin{pmatrix} \partial(\Im\,\pi^*) & \longto & \partial(\Geo(\Hc,T)) \\ b_{\tilde\gamma} & \longmapsto & b_\gamma\end{pmatrix} \]
   is well-defined and surjective.
  \item The map
\[ \pi_*\colon \begin{pmatrix} \partial^r(\Im\,\pi^*) & \longto & \partial^r(\Geo(\Hc,T)) \\ [b_{\tilde\gamma}]& \longmapsto & [b_\gamma]\end{pmatrix} \]
   is well-defined and surjective.
		\end{enumerate}
\end{corollary}
\begin{remark}
    The correct setting for Corollary \ref*{lemma:quotient} is that of graph covering. The map $\pi^*$ is the lift of paths. In particular, we could replace $(\Hc,d_T)$ by a Schreier graph.
\end{remark} 
\begin{proof} The results follow easily from Proposition \ref*{prop:switch}.
 \begin{enumerate}[leftmargin=8mm, label=(\arabic*)]
    \setcounter{enumi}{-1}
     \item We just have to check that $\tilde\gamma$ is a geodesic, which holds as $\pi$ is a submetry (i.e., $d_T(\pi(x),\pi(y))\le d_S(x,y)$) by construction of $T$.
     \item We only have to prove that $\pi_*$ is well-defined. Consider $\gamma,\eta\in\Geo(\Hc,T)$ and suppose that $b_{\tilde\gamma}=b_{\tilde\eta}$. For all $n\ge 0$, there exists $m\ge n$ such that
     \[ d_S(\tilde\gamma_m,\tilde\eta_n)=d_S(\tilde \eta_m,\tilde \gamma_n)=m-n \]
     hence $d_T(\gamma_m,\eta_n) = d_T(\pi(\tilde\gamma_m),\pi(\tilde\eta_n)) \le m-n$ as $\pi$ is a submetry. (The reverse inequality $d_T(\gamma_m,\eta_n)\ge m-n$ holds for any geodesics satisfying $\gamma_0=\eta_0$.) Using Proposition \ref{prop:switch}(1) once again, we conclude that $b_\gamma=b_\eta$ in $\partial(\Hc,T)$.
     \item Idem using Proposition \ref{prop:switch}(2) twice. \qedhere
 \end{enumerate}
\end{proof}
 
\subsection{Previous work on Abelian groups}

We recall several facts on horofunctions of abelian groups.
\begin{theorem}[{\cite[Prop.\ 3.5]{webster-winchester2006}}]\label{thm:websterwinchester}
Every horofunction on $\partial(\Z^r,S)$ is a Busemann point.
\end{theorem}

Therefore every horofunction on $\partial(\Z^r,d_S)$ is represented by an infinite geodesic word $\gamma=s_1s_2\cdots\in \Geo(\Z^r,S)$.
Let $\B=\Conv(S)\subset \R^r$ be the convex hull of $S$.
We state a lemma which is implicitly used in the previous work \cite{develin2002}.
\begin{definition}[\cite{develin2002,walsh2011}]
    Given an infinite word $\gamma\in S^\infty$, the set of \emph{directions} of $\gamma$ is the set of letters of $S$ with infinitely many occurrences in $\gamma$. We denote it by $D_\gamma$.
\end{definition}
\begin{lemma}\label{lemma:direction} Consider an abelian group $(\Z^r,S)$, and $\gamma\in S^\infty$ an infinite word.
\begin{enumerate}[leftmargin=8mm, label={\normalfont(\alph*)}]
    \item If $\gamma$ is a geodesic, the minimal face of $\B$ containing $D_\gamma$ is proper.. Equivalently, there exists a proper face $\F$ containing all but finitely many letters of $\gamma$.
    \item Reciprocally, if a proper face $\F\subset\B$ contains \emph{all} letters of $\gamma$, then $\gamma$ is geodesic.
\end{enumerate}
\end{lemma}
\noindent We will give its proof in more general setting (see Proposition \ref{proposition:faceinfiniteletters}). The minimal face containing $D_\gamma$ is the face \emph{associated to} $\gamma$. How these associated faces relate to Busemann points might not be clear, until the following result:

\begin{theorem}[{\cite[Section 4]{develin2002}}]\label{thm:abelian}
Let $\gamma_1,\gamma_2$ be infinite geodesic words in $(\Z^r,S)$.
\begin{enumerate}[leftmargin=8mm, label={\normalfont(\alph*)}]
    \item If $b_{\gamma_1}=b_{\gamma_2}$, then $\gamma_1$ and $\gamma_2$ are associated to the same face $\F$ of $\B$.
    \item Reciprocally, if \emph{all} letters of $\gamma_1$ and $\gamma_2$ belong to $\F$, then $b_{\gamma_1}=b_{\gamma_2}$.
    \item $b_{\gamma_1}, b_{\gamma_2}$ lie in the same orbit if and only if $\gamma_1,\gamma_2$ are associated to the same face $\F$. 
\end{enumerate}
\end{theorem}

\begin{remark}
    Part (c) follows from the first two parts and Lemma \ref{lemma:action}(b).
\end{remark}

\section{Horofunctions of $2$-step nilpotent groups}

\renewcommand{\H}{\mathcal{H}}
\newcommand{\dP}{d_\infty}

Throughout this section, we consider $\G$ a torsionfree $2$-step nilpotent group, and $S$ a finite symmetric generating set. We denote by $\Gbar$ the Malcev closure of $\G$. Note that $\Gbar$ can be identified with its Lie algebra $\mathfrak{g}\simeq \R^r$ via the exponential map. Hence the convex hull $\B=\Conv(S)\subset \Gbar$ is well-posed. Let $\Pr\colon \Gbar \to \Gbar/[\Gbar,\Gbar] \simeq \R^d$ be the abelianization map. Since $\Pr$ is a homomorphism, $\Pr(\B)$ is the convex hull of $\Pr(S)$.
\subsection{Infinite geodesic words}

The main goal of this subsection is to extend Lemma \ref{lemma:direction} in the setting of finitely generated, torsionfree, $2$-step nilpotent groups. Using the previous notations, we prove

\begin{proposition}\label{proposition:faceinfiniteletters}
Let $\G$ be a torsionfree $2$-step nilpotent group, and $S$ a finite symmetric generating set. 
For any infinite geodesic word $\gamma\in \Geo(\G,S)$, there exists a proper face $\F\subset \Pr(\B)$ such that all but finitely many letters of $\gamma$ projects to $\F$.
\end{proposition}

\newpage
    
To prove this proposition, we need to recall the work of Stoll in \cite{Stoll_metric}.

\begin{definition}
An $\R$-word is an expression of the form
\[ w = s_1^{\lambda_1} \cdot s_2^{\lambda_2}\cdot\ldots\cdot s_k^{\lambda_k} \]
with $s_i\in S$, $s_{i+1}\ne s_i$ and $\lambda_i>0$. For each $\R$-word $w$, we define the followings.
\begin{itemize}[leftmargin=6mm]
    \item Its \emph{length}: $\abs{w}_S = \sum_{i=1}^k \lambda_i$,
    \item The \emph{total exponent} of a letter $s\in S$ in a word $w$: $\abs{w}_s = \sum_{i:s_i=s}\lambda_i$,
    \item Its \emph{coarse length}:  $k(w)=k$.
\end{itemize}

The set of $\R$-words is denoted by $S_\R^*$. Note that each $\R$-word $w\in S_\R^*$ represents a well-defined element $\bar w\in \Gbar$. For each $g\in\Gbar$, we define
\[ \norm{g}_\Stoll = \inf\left\{\abs w_S \,\big|\, w\in S_\R^* \text{ and } \bar w=g \right\}.\]
\end{definition}

\begin{remark}
The Stoll metric is the Lie group analogue of the word metric:
\[\norm{g}_S \coloneqq \inf\left\{\abs w_S \,\big|\, w\in S^* \text{ and } \bar w=g \right\}=d(e,g).\]
Note that the inequality $\norm{g}_\Stoll\le \norm{g}_S$ holds for any $g\in \G$.
\end{remark}

\medskip

\begin{lemma}[{\cite[Lemma 3.3]{Stoll_metric}}] \label{sec4:lem_coarse}
There exists a constant $K$ such that, for each $\R$-word $w\in S_\R^*$, there exists another $w'\in S_\R^*$ satisfying the following conditions:
\begin{itemize}[leftmargin=6mm]
    \item Both words represent the same elements $\overline w=\overline{w'}$ in $\Gbar$.
    \item For each letter $s\in S$, we have $\abs{w}_s\ge \abs{w'}_s$. In particular, $\abs{w}_S\ge \abs{w'}_S$.
    \item $w'$ has uniformly bounded coarse length, precisely $k(w')\le K$.
\end{itemize}
If $w$ is geodesic, then so is $w'$, and we have $\abs{w}_s=\abs{w'}_s$ for all $s\in S$.
\end{lemma}

\medskip

\begin{theorem}[{\cite[Theorem 4.5]{Stoll_metric}}] \label{sec4:rough}
There exists $C\ge 0$ such that
\[ \forall g\in \G,\quad \norm{g}_\Stoll \le \norm{g}_S \le \norm{g}_\Stoll + C. \]
\end{theorem}

\bigskip

We have laid down everything needed to prove the main proposition of this section.

\begin{proof}[Proof of Proposition \ref{proposition:faceinfiniteletters}]For every pair $\{s,t\}\subset S$ such that $\Pr(s)$ and $\Pr(t)$ do \textbf{not} lie on a common face of $\Pr(\B)$, there exists an $\R$-word $v=v(s,t)$ such that $\Pr(st)=\Pr(\bar v)$ and $\abs{v}_S < 2 = \abs{st}_S$. We denote by $z(s,t)\in [\Gbar,\Gbar]$ the element satisfying $st = \bar v z$. Let $\delta>0$ be the minimum value of $2 - \abs{v(s,t)}_S$ over all such pairs.

\bigbreak

Consider $\gamma\in S^\infty$ and let us partition $S=D_\gamma\sqcup D_\gamma^c$. We suppose that no proper face of $\Pr(\B)$ contains $\Pr(D_\gamma)$, and provide a finite subword $v$ of $\gamma$ such that 
\[ \norm{\bar v}_\Stoll<\abs v_S-C, \]
hence neither $v$ nor $\gamma$ can be geodesic w.r.t.\ the word metric (by Theorem \ref*{sec4:rough}). 

\bigbreak

By hypothesis, $\gamma$ can be decomposed as
\[ \gamma = w_0 \cdot s_1u_1t_1 \cdot w_1\cdot s_2u_2t_2 \cdot w_2 \cdot \ldots \cdot w_{n-1} \cdot s_nu_nt_n \cdot w_n \cdot w_\infty \]
where
\begin{itemize}[leftmargin=8mm]
    \item $u_i,w_i\in S^*$ are finite words, with all the (finitely many) instances of $D_\gamma^c$ in $\gamma$ appearing in $w_0$, and $w_\infty\in S^\infty$ is an infinite word,
    \item $s_i$, $t_i\in D_\gamma$ are individual letters such that the projections $\Pr(s_i),\Pr(t_i)$ do \textbf{not} lie on a common face for each $i\in\{1,2,\ldots,n\}$,
    \item for each $s\in D_\gamma$, we have $\abs{w_n}_s \ge K\abs{u_1u_2\ldots u_n}_s + C$,
    \item the parameter $n$ can be made arbitrarily large. 
\end{itemize}

\bigbreak

\noindent The subword we are looking for is
\[ v= s_1u_1t_1 \cdot w_1\cdot s_2u_2t_2 \cdot w_2 \cdot \ldots \cdot w_{n-1} \cdot s_nu_nt_n \cdot w_n. \]
As $\Gbar$ is $2$-step nilpotent, the derived subgroup $[\Gbar,\Gbar]$ is central, we can rewrite
    \begin{align*}
        \overline v & = s_1u_1t_1w_1 \cdot\ldots\cdot s_nu_nt_nw_n \\
        & = [s_1,u_1]u_1s_1t_1 w_1  \cdot\ldots\cdot [s_n,u_n]u_ns_nt_nw_n \\
        & = u_1v_1w_1 \cdot \ldots \cdot u_nv_n \cdot w_n \cdot [s_1,u_1]\ldots [s_n,u_n] \cdot z_1\ldots z_n,
    \end{align*}
    where $v_i=v(s_i,t_i)$ and $z_i=z(s_i,t_i)$ defined earlier (the key properties being $s_it_i=v_iz_i$ and $\norm{v_i}_\Stoll\le 2-\delta$). In particular, we have
    \[ \norm{u_1v_1w_1 \cdot \ldots \cdot u_nv_n}_\Stoll \le \abs{s_1u_1t_1w_1 \cdot\ldots\cdot s_nu_nt_n}_S - n\delta, \]
    and
    \[ \norm{z_1\ldots z_n}_\Stoll = O(\sqrt n),\]
    so we only need good estimates on $\norm{w_n \cdot [s_1,u_1]\ldots [s_n,u_n] }_\Stoll$, i.e., a short $\R$-word representing $w_n \cdot [s_1,u_1]\ldots [s_n,u_n]$. Using Lemma \ref{sec4:lem_coarse}, we can find $w'_n\in S_\R^*$ such that $\bar{w_n}=\bar w'_n$ and $k(w'_n)\le K$. Two cases present themselves:
\begin{itemize}[leftmargin=6mm]
    \item If $\abs{w'_n}_s<\abs{w_n}_s-C$ for some $s\in S$, then $\abs{w'_n}_S<\abs{w_n}_S-C$ too and
    \[ \norm{v}_\Stoll \le \abs{s_1u_1t_1 \cdot w_1\cdot s_2u_2t_2 \cdot w_2 \cdot \ldots \cdot w_{n-1} \cdot s_nu_nt_n \cdot w'_n}_S<\abs{v}_S-C, \]
    so can conclude without much of the conversation above.
    
    \item Otherwise, we have $\abs{w'_n}_s\ge \abs{w_n}_s-C\ge K\abs{u_1u_2 \ldots u_n}_s$ for all $s$. Knowing that $k(w'_n)\le K$, the pigeonhole principle implies that the $\R$-word $w'_n$ contains a power $s^{\lambda_s}$ with $\lambda_s\ge \abs{u_1u_2\ldots u_n}_s$ for each $s\in D_\gamma$.
    
    From $w'_n$, we construct a short word $w''_n$ representing $w_n \cdot [s_1,u_1]\ldots [s_n,u_n]$. For each generator $s\in S_\infty$, we replace in $w'_n$ the previous power $s^{\lambda_s}$ by the $\R$-word
    \[ s_n^{-\abs{u_n}_s/\lambda_s}\ldots s_1^{-\abs{u_1}_s/\lambda_s}\cdot s^{\lambda_s} \cdot  s_1^{\abs{u_1}_s/\lambda_s} \ldots s_n^{\abs{u_n}_s/\lambda_s}\]
    so that 
    $\abs{w''_n}_S = \abs{w'_n}_S + 2 \sum_{s\in S_\infty} \frac1{\lambda_s} \abs{u_1u_2\ldots u_n}_s \le \abs{w_n}_S + 2\abs{D_\gamma}$. Putting everything together, the triangle inequality gives
    \begin{align*}
        \norm{\bar v}_\Stoll
        & \le \norm{u_1v_1w_1 \cdot \ldots \cdot u_nv_n}_\Stoll + \norm{w_n \cdot [s_1,u_1]\ldots [s_n,u_n] }_\Stoll + \norm{z_1\ldots z_n}_\Stoll \\
        & \le \Big(\abs{s_1u_1t_1w_1 \cdot\ldots\cdot s_nu_nt_n}_S - n\delta\Big) + \Big(\abs{w_n}_S + 2\abs{D_\gamma}\Big) + O(\sqrt n) \\
        & = \abs{v}_S - n\delta + O(\sqrt n) + 2\abs{D_\gamma} \\
        & < \abs{v}_S-C
    \end{align*}
    for $n$ large enough. \qedhere
\end{itemize}
\end{proof}

\subsection{Orbits of Busemann points}

In this section, we consider $\Hc$ a finitely generated, torsionfree, $2$-step nilpotent group, with the additional condition $[\Hc,\Hc]=\la z\ra\simeq\Z$. These are lattices inside $\H^k_\R\times\R^\ell$ (see e.g. \cite[Lem.\ 7.1]{Stoll-rational}), where $\H^k_\R$ is the $(2k+1)$-dimensional Heisenberg Lie groups
\[
\Hc^k_\R = \left\{
\begin{pmatrix}
1 & \mathbf a & c \\
& \mathrm{Id}_k & \mathbf b^t \\
& & 1
\end{pmatrix}
\;\middle|\; \mathbf a,\mathbf b\in \R^k,\, c\in\R \right\}.
\]
The goal is to generalize Theorem \ref{thm:abelian} and classify orbits of Busemann points in $\partial(\Hc,d_S)$, for any finite symmetric generating set $S$. In view of Proposition \ref{proposition:faceinfiniteletters}, we might expect that orbits of Busemann points are once again classified by proper faces of $\Pr(\B)$. However, the conclusion is not so straightforward, and we often need more information to chose in which orbit a given Busemann point $b_\gamma$ fits. Namely, for each geodesic ray $\gamma\in\Geo(\Hc,S)$, we define
\begin{itemize}[leftmargin=8mm]
    \item $D_\gamma\subseteq S$, the set of letters that appears infinitely often in $\gamma$.
    \item $\E_\gamma\subseteq \B$, the minimal face of $\B$ containing $D_\gamma$.
    \item $\F_\gamma\subseteq \Pr(\B)$, the minimal face of $\Pr(\B)$ containing $\Pr(D_\gamma)$. 
\end{itemize}
Note that $\E_\gamma$ is a face of $\B\cap\Pr^{-1}(\F_\gamma)$, and $\F_\gamma$ is a \emph{proper} face of $\Pr(\B)$.
\begin{definition}
A face $\F\subset \Pr(\B)$ is \emph{commutative} if $\Ac_\R=\la \Pr^{-1}(\F)\ra\le \Hc_\R$ is abelian.
\end{definition}

\begin{theorem} \label{theorem:nheisenberg orbits}
Let $(\Hc,S)$ be a $2$-step nilpotent group with $[\Hc,\Hc]=\la z\ra\simeq\Z$. Consider two geodesic rays $\gamma_1,\gamma_2\in\Geo(\Hc,S)$.
\begin{enumerate}[leftmargin=8mm, label={\normalfont(\alph*)}]
    \item If $b_{\gamma_1}=b_{\gamma_2}$, then $\F_{\gamma_1}=\F_{\gamma_2}$. Moreover, if this face is commutative, then $\E_{\gamma_1}=\E_{\gamma_2}$.
    \item Reciprocally, if
    \begin{enumerate}[leftmargin=6mm, label={\normalfont(\arabic*)}]
        \item all the letters of $\gamma_1,\gamma_2$ belong to $\E$ with $\E$ commutative, or
        \item all the letters of $\gamma_1,\gamma_2$ project to $\F$ with $\F$ non-commutative,
    \end{enumerate}
    then $b_{\gamma_1}=b_{\gamma_2}$.
    \item The Busemann points $b_{\gamma_1},b_{\gamma_2}$ are in the same orbit if and only if
    \begin{enumerate}[leftmargin=6mm, label={\normalfont(\arabic*)}]
        \item $\F_{\gamma_1}=\F_{\gamma_2}$ is commutative and $\E_{\gamma_1}=\E_{\gamma_2}$, or
        \item $\F_{\gamma_1}=\F_{\gamma_2}$ is non-commutative.
    \end{enumerate}
\end{enumerate}
In particular, $\partial(\Hc,d_S)$ contains only finitely many orbits of Busemann points {\normalfont(}since $\B$ has only finitely many faces{\normalfont)}, hence countably many Busemann points.
\end{theorem}

\bigskip

Before coming to a proof of Theorem \ref*{theorem:nheisenberg orbits}, let us do some anagrams!
\begin{definition}
    For $w\in S^*$, we define
    \[ Z_w = \{ g\in\la z\ra : \exists w'\text{ a reodering of } w\text{ such that }\bar w'=\bar wg \} \subseteq [\Hc,\Hc]\simeq \Z.\]
    This definition extends to infinite words $\gamma\in S^\infty$, letting $Z_\gamma = \bigcup_{n\ge 0} Z_{\gamma_n}$, where $\gamma_n\in S^*$ is the prefix of length $n$ of $\gamma$.
\end{definition}
\begin{lemma} \label{lem:non_extremal}
    Let $(\Hc,S)$ be a $2$-step nilpotent group with $[\Hc,\Hc]=\la z\ra\simeq\Z$. Consider $s,t\in S$ such that $[s,t]\ne e$. Then, for any word $w\in S^*$ containing the letters $s,t,s,t$ in that order, the set $Z_w\subset \Z$ contains both positive and negative numbers.
\end{lemma}
\begin{proof}
We will show the existence of $w'$ such that $\bar w'=\bar wz^a$ with $a>0$. Without lost of generality, we suppose $w=s\,v_1\,t\,v_2\,s\,v_3\,t$. We split into different cases:
    \begin{itemize}[leftmargin=8mm]
        \item If $[s,v_1tv_2]\ne e$, then either $w'=s^2v_1tv_2v_3t$ 
        or $w'=v_1tv_2s^2v_3t$ satisfies. (More precisely, the first satisfies if $[s,v_1tv_2]= z^a$ with $a>0$, and the second otherwise.)
        \item If $[t,v_2sv_3]\ne e$, then either $w'=sv_1t^2v_2sv_3$ 
        or $w'=sv_1v_2sv_3t^2$ satisfies.
        \item Otherwise, we have $\bar w=sv_1t^2v_2sv_3$. Moreover
        \[ [s,v_1t^2v_2] = [s,v_1tv_2] \cdot [s,t] \ne e, \]
        hence either $w'=s^2v_1t^2v_2v_3$ or $w'=v_1t^2v_2s^2v_3$ satisfies. \qedhere
    \end{itemize}
\end{proof}

\begin{lemma} \label{lem:interval}
    Let $(\Hc,S)$ be a $2$-step nilpotent group with $[\Hc,\Hc]=\la z\ra\simeq\Z$. Consider an infinite word $\gamma\in S^\infty$ satisfying $\gamma\in D_\gamma^\infty$. Then $Z_\gamma = \la [s,t] : s,t\in D_\gamma\ra$.
\end{lemma}

\begin{proof} The inclusion $Z_\gamma \subseteq \la [s,t] : s,t\in D_\gamma\ra$ is clear. In particular, the statement is true if $[s,t]=e$ for all $s,t\in D_\gamma$. From now on, we suppose the existence of $s,t\in D_\gamma$ such that $[s,t]\ne e$. Let's make some observations on $Z_\gamma$:
\begin{itemize}[leftmargin=8mm]
    \item Lemma \ref{lem:non_extremal} implies that $Z_\gamma$ (seen as a subset of $\Z$) is unbounded above and below.
    \item Consider $\delta=\max\{ \abs a : \exists s,t\in S,\, [s,t]=z^a\}$. Then $Z_\gamma\subseteq\Z$ is a $\delta$-net, meaning thickening $Z_\gamma$ gets you $Z_\gamma+[-\delta,\delta]=\R$. Indeed, every reordering of $\gamma_n$ can be obtained by successively permuting pairs of consecutive letters $s,t$, hence changing the value by a commutator $[s,t]=z^a$ with $\abs a\le \delta$ at each step.
    \item There exists $u\in D_\gamma^*$ such that $Z_u\supseteq \la [s,t] : s,t\in D_\gamma\ra\cap [-\delta,\delta]$. For instance, if one enumerates $D_\gamma=\{s_1,s_2,\ldots,s_n\}$, then
    \[ u = (s_1s_2)^N(s_1s_3)^M\ldots (s_1s_n)^M(s_2s_3)^M\ldots (s_{n-1}s_n)^M\]
    satisfies when $M$ is large enough.
\end{itemize}
Let $\gamma=\gamma_n\omega_n$ with $\omega_n$ the infinite suffix. For $n$ large enough, we can reorder $\gamma_n$ into a word starting with $u$, say $\gamma'_n=uv$. Fix $z^a\in\la [s,t]:s,t\in D_\gamma\ra$ s.t.\  $\bar\gamma'_n=\bar\gamma_nz^a$. Then
\[ Z_\gamma = a + Z_{\gamma_n'\omega_n} \supseteq a + Z_u + Z_{v\omega_n} = \la [s,t]:s,t\in D_\gamma\ra \]
as $Z_u\supseteq \la [s,t] : s,t\in D_\gamma\ra\cap [-\delta,\delta]$ and $Z_{v\omega_n}$ is a $\delta$-net.
\end{proof}

\bigskip

\begin{proof}[Proof of Theorem \ref{theorem:nheisenberg orbits}]

\medbreak

\noindent (a) We suppose that $b_{\gamma_1}=b_{\gamma_2}$.

\bigbreak

\noindent$\blacktriangleright$ By Proposition \ref{prop:switch},
there is a geodesic word $\gamma$ having infinitely many intersections with both $\gamma_1$ and $\gamma_2$. We consider the projections $\Pr(\gamma_1),\Pr(\gamma)\in \Pr(S)^\infty$ in the abelianization. By Proposition \ref{proposition:faceinfiniteletters}, all but finitely many letters of $\gamma$ project to the face $\F_{\gamma}\subset \partial\B$, in particular $\Pr(\gamma)$ is eventually geodesic.
Since $\gamma_1$ and $\gamma$ have infinitely many intersections,
their projections $\Pr(\gamma_1),\Pr(\gamma)$ also have infinite intersection. In particular, $\Pr(\gamma_1)$ is also eventually geodesic.
By Proposition \ref{prop:switch} and Theorem \ref{thm:abelian},
we conclude $\F_{\gamma_1}=\F_{\gamma}$. We prove $\F_{\gamma_2}=\F_{\gamma}$ in the same way, hence $\F_{\gamma_1}=\F_{\gamma_2}$.

\bigbreak

\noindent$\blacktriangleright$ Suppose that $\F=\F_{\gamma_1}=\F_{\gamma_2}=\F_\gamma$ is commutative.

\medbreak

Let $S_\F\coloneqq S\cap\Pr^{-1}(\F)$. Then $S_\F^{\pm 1}$ generates a discrete abelian subgroup $\mathcal{A}$ in $\Hc$. Denote by $\mathcal{A}_\R$ the Mal'cev completion of $\mathcal{A}$, and $\Conv_{\mathcal{A}}(S_\F^{\pm 1})$ the convex hull of $S_\F^{\pm 1}$ in $\mathcal{A}_\R$. By Theorem \ref{thm:abelian},
 the orbits of Busemann points of $(\mathcal{A},S_\F^{\pm 1})$ are in one-to-one correspondence with the faces of $\Conv_{\mathcal{A}}(S_\F^{\pm 1})$.
    Note that the faces of $\Conv_{\mathcal{A}}(S_\F^{\pm 1})\cap\Pr^{-1}(\F)$ coincides with the faces of $\B\cap\Pr^{-1}(\F)$.

\medbreak

As we have just seen, all the letters of $\gamma_1,\gamma_2,\gamma$ belong to $S_\F$ after some time $N$. In particular, these paths are geodesic rays in $(\mathcal{A},S_\F^\pm)$ after the time $N$. As $\gamma$ still intersects $\gamma_1$ and $\gamma_2$ after time $N$, we have $b_{\gamma_1}=b_{\gamma_2}$ in $\partial(\mathcal{A},S_\F^\pm)$ by Proposition \ref{prop:switch}, hence $\E_{\gamma_1}=\E_{\gamma_2}$ by Theorem \ref{thm:abelian}.

\medbreak

\noindent (b1) Let $\F$ be a commutative face and $\E$ be a face of 
\[\Conv_{\mathcal{A}}(S_\F^{\pm 1})\cap {\Pr}^{-1}(\F)=\B\cap {\Pr}^{-1}(\F). \]
Consider infinite words $\gamma_1,\gamma_2\in S^\infty$ such that all the letters belong to $\E$. In particular, we have $\gamma_1,\gamma_2\in\Geo(\mathcal{A},S_\F^\pm)$. By Theorem \ref{thm:abelian}, $\gamma_1$ and $\gamma_2$ yield the same Busemann points in $(\mathcal{A},S_\F^\pm)$. By Proposition \ref{prop:switch},
there exists $\gamma\in\Geo(\mathcal{A},S_\F^\pm)$ intersecting with $\gamma_1, \gamma_2$ infinitely many times. Necessarily $\gamma\in S_F^\infty$. Since $\Pr(S_\F)\subset \F$, $\gamma_1, \gamma_2$ and $\gamma$ are also geodesic rays in $(\H,S)$. By Proposition \ref{prop:switch} again, we have $b_{\gamma_1}= b_{\gamma_2}$ in $\partial(\Hc,d_S)$.

\medbreak

\noindent(b2) Consider $\gamma,\eta\in \Geo(\H,S)$ such that $\F=\F_\gamma=\F_\eta$ is a non-commutative face and \emph{all} letters of $\gamma,\eta$ belong to $S_\F$. We prove that $b_\gamma=b_\eta$. More specifically, we fix a prefix $\eta_m$ of $\eta$, and attempt to extend it into a longer geodesic joining the path $\gamma$. Fix $n$ such that all occurrences of letters of $S_\F\setminus D_\gamma$ in $\gamma=\gamma_n\omega_n$ belong to $\gamma_n$.

\medbreak

Observe that $\la D_\gamma \ra$ is a \emph{finite-index} subgroup of $\la S_\F\ra$. Indeed, we have
\begin{center}
    \begin{tikzcd}
        1 \arrow[r] & \la D_\gamma\ra\cap \la z\ra \arrow[r,hook] \arrow[d,hook,"(1)"] & \la D_\gamma\ra \arrow[r,two heads, "\Pr"] \arrow[d,hook] & \Pr\la D_\gamma\ra \arrow[r] \arrow[d,hook,"(2)"] & 1 \\
        1 \arrow[r] & \la S_\F\ra\cap \la z\ra \arrow[r,hook] & \la S_\F \ra \arrow[r,two heads, "\Pr"] & \Pr\la S_\F\ra \arrow[r] & 1
    \end{tikzcd}
\end{center}
where the inclusion (1) has finite-index as both are non-trivial subgroup of $\la z\ra$ (since the face $\F$ is non-commutative), and (2) has finite-index as both are lattices in the hyperplane of $\Hc_\R/[\Hc_\R,\Hc_\R]$ supporting $\F$. Therefore, there exists $u_1\in S_\F^*$ such that $\bar \gamma_n^{-1}\cdot \overline{\eta_m u_1}\in \la D_\gamma\ra$. Using \cite[Lemma 4.2]{walsh2011}, there exist $u_2,v\in D_\gamma^*$ such that
\[ \overline{\eta_m u_1\cdot u_2} = \overline{\gamma_n\cdot v}. \]

As every letter of $D_\gamma$  appears infinitely many times in $\omega_n$, there exists a prefix $w_1\in D_\gamma^*$ of $\omega_n$ (so $\gamma_p=\gamma_nw_1$ for $p>n$) containing all the letters of $v$. Then $w_1$ can be reordered as $w_1'=vu_3$ for some $u_3\in D_\gamma^*$, and $\bar w_1'=\bar w_1g$ with $g\in \la [s,t]:s,t\in D_\gamma\ra$.

\medbreak

Using Lemma \ref{lem:interval}, we have $g^{-1}\in Z_{\omega_p}$, hence there exists a prefix $w_2\in D_\gamma^*$ of $\omega_p$ (so $\gamma_q=\gamma_pw_2$ for $q>p$) which can be reordered as $w_2'$ with $\bar w_2'=\bar w_2g^{-1}$.

\medbreak

To summarize, we have $\eta_mu_1u_2u_3w_2'\in S_\F^*$ a geodesic extension of $\eta_m$ such that
\[ \overline{\eta_mu_1u_2u_3w_2'} = \overline{\gamma_nvu_3w_2'} = \overline{\gamma_nw_1'w_2'} = \overline{\gamma_nw_1w_2} = \bar\gamma_q. \]
This can be repeated, extending $\gamma_q$ to join $\eta$, hence $b_\gamma=b_\eta$ by Proposition \ref{prop:switch}.

\bigbreak

\noindent (c) follows from parts (a) and (b), together with Lemma \ref{lemma:action}(b).
\end{proof}

In order to generalize these results to general finitely generated $2$-step nilpotent groups, a first tool could be Lemma \ref{lemma:quotient}, as we could reduce the problem to $(N_{2,r},S_r)$, the free $2$-step nilpotent group of $r$ with its standard generating set. Indeed, for any $2$-step nilpotent group $\G$ generated by $S$, there exists a surjective homomorphism $f\colon N_{2,r}\to\G$ such that $f(S_r)=S$, for $r=\abs S$. We isolate the following conjecture:
\begin{conjecture}
    The set of Busemann points of the free nilpotent group of step $2$ and rank $r$, with respect to the standard generating set, is countable.
    \end{conjecture}

\subsection{Non-Busemann points}

As we saw in Theorem \ref{thm:websterwinchester}, all horofunctions of an abelian group are Busemann points. In this section, we shall see that this fails for $2$-step nilpotent groups. For the simplicity, we focus on lattices $\Hc$ inside the $3$-dimensional Heisenberg Lie group $\Hc^1_\R$. We identify $\H_\R^1$ with its Lie algebra $V_1\oplus V_2$, with $V_1\simeq \R^2$ and $V_2\simeq \R$. The operation is then
\[
(a,b,c)\cdot (a',b',c') = \left( a+a',\; b+b',\; c+c'+\frac12(ab'-a'b)\right).
\]
Moreover, we identify $\H_\R^1/[\H_\R^1,\H_\R^1]\simeq\R^2$ with $V_1\subset \H_\R^1$. Let $\Pr(S)\subset V_1$ be the projection of a finite generating set $S$. Note that $\Pr(S)$ is a Lie generating set of $\H_\R^1$. We shall denote by $\|\cdot\|_\infty$ be the Stoll distance on $\H_\R^1$ associated to $\Pr(S)$. Note that, since $\Pr(S)$ is in the first layer of the Heisenberg group,
the induced distance is a left invariant sub-Finsler distance admits dilations,
that is, the equality $\|(ta,tb,t^2c)\|_\infty=t\|(a,b,c)\|_\infty$ holds for all $(a,b,c)\in \H_\R^1$ (cf.\ \cite[Theorem 2]{berestovskii1988}). Although the distance $\|\cdot\|_\infty$ is different from the Stoll distance $\|\cdot\|_{\Stoll}$,
it still approximates the word metric:
\begin{theorem}[{\cite[Theorem 0.1]{krat2002}}]
    There is a constant $C>0$ such that
    \[\forall g\in\Hc,\quad \|g\|_\infty-C\le \|g\|_S\le \|g\|_\infty+C.\]
\end{theorem}

On the other hand, in \cite{fisher2021sub}, Fisher and Nicolussi-Golo characterize all the horofunctions of $(\H_\R^1,\|\cdot\|_\infty)$.
We re-state their result up to bounded functions.

\begin{itemize}[leftmargin=8mm]
\item $\Pr(\B)\coloneqq \Conv(\Pr(S))$,
    \item $\mathsf v_1,\ldots, \mathsf v_{2N}$ the vertices of $\B$ numbered counterclockwise. By convention $\mathsf v_0=\mathsf v_{2N}$.
    \item $e_k \coloneqq \mathsf v_k-\mathsf v_{k-1}$.
    \item $\omega\colon V_1\times V_1\to\R$ defined by $\omega\big((a,b),(a',b')\big)=a'b-ab'$.
\item $\alpha_k(a,b)=\frac{\omega(e_k,v)}{\omega(e_k,\mathsf v_k)}$ for $v\in V_1$.
\end{itemize}
\begin{theorem}[{\cite[Theorem 5.2]{fisher2021sub}}]
        The horofunctions of a polygonal sub-Finsler Heisenberg group $(\H_\R^1,d_{sF})$ are, up to bounded functions, classified as follows:
        \begin{description}
            \item[\quad Vertical] $\varphi(v,c)=-\|v\|_{\Pr(\B)}\coloneqq\min\{\lambda\ge 0:x\in\lambda\Pr(\B)\}$,
            \item[\quad Non-vertical] $\varphi(v,c)=r\alpha_k(v)+(1-r)\alpha_{k-1}(v)$ for some $r\in[0,1]$,
            \item[\quad Mixed]
            \[\varphi(v,c)=\begin{cases}
                \alpha_i(v) & \text{if}~~\omega(\mathsf v_i,v)\le0 ~~ (\text{resp. }\ge0),\\
                r\alpha_i(v)+(1-r)\alpha_{i-1}(v) & \text{if}~~\omega(\mathsf v_i,v)\ge 0 ~~ (\text{resp. }\le 0),
            \end{cases}\]
            or
            \[\varphi(v,c)=\begin{cases}
                \alpha_{i-1}(v) & \text{if}~~ \omega(\mathsf v_i,v)\le0 ~~ (\text{resp. }\ge0),\\
                r\alpha_i(v)+(1-r)\alpha_{i-1}(v) & \text{if}~~ \omega(\mathsf v_i,v)\ge0 ~~ (\text{resp. }\le0),
            \end{cases}\]
            for some $r\in[0,1]$.
        \end{description}
    \end{theorem}
Combined with Proposition \ref{proposition:equivalencehorofunction},
we obtain the following corollary.
  
\begin{corollary}\label{cor:Heisenberg}
    For a lattice of $3$-dimensional Heisenberg group $\H<\H_\R^1$ with a generating set $S$,
    its Cayley graph $(\H,S)$ has uncountably many horofunctions.
\end{corollary}

It can be expected that the horofunction boundary of any $2$-step nilpotent group is uncountable, unless the group is virtually abelian. It is known that a Cayley graph is roughly isometric to a Carnot group with a homogeneous sub-Finsler metric if its Mal'cev completion is ideal (see \cite{Tashiro2022}),
and the horofunctions of a Carnot group is characterized by the Pansu derivative of the distance function (see \cite[Lemma 2.3]{fisher2021sub}). 
However, if the Mal'cev completion is not ideal,
then there are examples of Cayley graphs which are not roughly isometric to any Carnot group (see \cite{Breuillard-Ledonne2013}).
In this case, horofunctions may not be given by the Pansu derivative of the distance function.
To describe $\partial^r(\G,d_S)$ for general $2$-step nilpotent groups, we need a characterization of horofunctions which does not use Pansu derivatives.

\section{The Cartan group}

\subsection{A model for the Cartan group}

We recall a model for the Cartan group $\Cc$. We first define a nilpotent Lie group $\Cc_\R$, in which $\Cc$ sits as a cocompact lattice. Elements of $\Cc_\R$ can be understood geometrically: They are equivalence classes of absolutely continuous paths in $\R^2$ starting from $(0,0)$. For any path $g$, we define three parameters:
\begin{itemize}[leftmargin=8mm, label=\textbullet]
	\item[(1)] its second endpoint $\hat g=(x_g,y_g)\in\R^2$.
	\item[(\textbullet)] a distribution of winding numbers. First, we get a closed path $g_c$ by concatenating $g$ with the segment back from $\hat g$ to $(0,0)$. Then the function $w_g\colon\R^2\setminus \Im(g_c)\to \Z$ is defined as $w_g(x,y)=$ the winding number of $g_c$ around $(x,y)$.
	\item[(2)] its total algebraic (or signed) area \vspace*{1mm}
	\[ A(g) = \displaystyle \iint_{\R^2} w_g(x,y)\;\dif x\,\dif y \in \R. \vspace*{-1mm} \]
	\item[(3)] its (non-normalized) barycenter (or center of gravity)
	\[ B(g) = \displaystyle\iint_{\R^2} (x,y) \cdot w_g(x,y) \;\dif x\,\dif y \in \R^2. \]
\end{itemize}
Two paths $g,h$ are equivalent if they share same endpoint $\hat g=\hat h$, same algebraic area $A(g)=A(h)$ and same \say{barycenter} $B(g)=B(h)$.

\medbreak

\begin{proposition}[{\cite[Proposition 4.1]{Intermediate}}]Given two paths $g,h$ in $\R^2$, their concatenation $gh$ has parameters
	\begin{align*}
	\widehat{gh} \hspace*{1.5mm} & = \widehat g + \widehat h \\
	A(gh) & = A(g) + A(h) + \frac12\det\!\big(\hat g,\hat h\big) \\
	B(gh) & = B(g)+ B(h) + \hat g \cdot A(h) + \frac13(2\hat g+ \hat h)\cdot\frac12\det\!\big(\hat g,\hat h\big)
	\end{align*}
	As a corollary, the operation \say{concatenation} passes to the quotient. With the empty path as neutral element and reverse paths as inverses, this defines a group.
\end{proposition}

\bigbreak

We consider the lattice $\Cc=\la x,y\ra$, where $x$ (resp.\ $y$) is the straight segment from $(0,0)$ to $(1,0)$ (resp.\ to $(0,1)$). This group coincides with $N_{3,2}$ the free $3$-step nilpotent group of rank $2$. Indeed, our group $\Cc$ is $3$-step nilpotent and has rank $2$, so $\Cc$ is a quotient of $N_{3,2}$. Moreover, for any proper quotient $G$ of $N_{3,2}$, the subgroup $[G,[G,G]]$ is a proper quotient of $\Z^2$. In contrast $[\Cc,[\Cc,\Cc]]$ contains a copy of $\Z^2$, so $\Cc\simeq N_{3,2}$.

\subsection{Construction}

\renewcommand{\SS}{\mathbb{S}}
\newcommand{\atob}{[a\!:\!b]}

We are going to define a Busemann point $\gamma_u$ for each direction $u\in\SS^1$ in the plane. Directions can also be parametrized by equivalence classes $\atob\in \big(\R^2\setminus\{(0,0)\}\big)\big/\,\R^+$, we will sometimes use this second point of view.

\medbreak

Note that, in order to define an infinite path $\gamma=(\gamma_n)_{n\ge 0}$, we only need to fix $\gamma_0\in\Cc$, and its projection $\hat\gamma=(\hat \gamma_n)_{n\ge 0}\in\Z^2\simeq\Cc/[\Cc,\Cc]$. The curve $\gamma$ is the unique lift. In this case, we take $\gamma_{u,0}=e$, and $(\hat\gamma_{u,n})_{n\ge 0}$ which best approximates the ray $\R^+u$. Precisely, unit squares in the grid fall into two categories, depending whether the ray passes above or below their center. $\hat\gamma_u$ is the boundary of the two regions formed. 
	
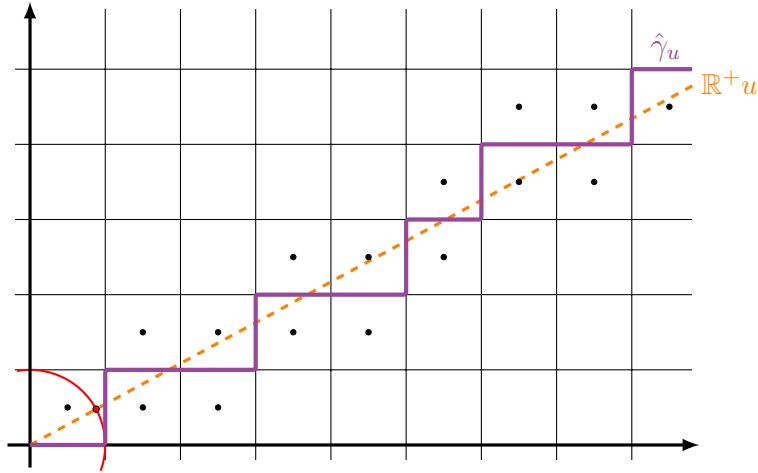
\begin{center}
		\begin{tikzpicture}[scale=1]
		\draw[thin] (-.2,-.2) grid (8.8,5.8);
		
		\draw[red, thick, domain=-20:100] plot ({cos(\x)}, {sin(\x)});
		
		\draw[-latex, very thick] (-.3,0) -- (8.9,0);
		\draw[-latex, very thick] (0,-.3) -- (0,5.9);
		
		\draw[orange, dashed, very thick] (0,0) -- (8.9,4.83);
		
		\draw[fill=red] ({8.9/sqrt(8.9^2+4.83^2)},{4.83/(sqrt(8.9^2+4.83^2))}) circle (1.2pt);
		
		\draw[fill=black] (.5,.5) circle (1pt);
		\draw[fill=black] (1.5,.5) circle (1pt);
		\draw[fill=black] (1.5,1.5) circle (1pt);
		\draw[fill=black] (2.5,0.5) circle (1pt);
		\draw[fill=black] (2.5,1.5) circle (1pt);
		\draw[fill=black] (3.5,1.5) circle (1pt);
		\draw[fill=black] (3.5,2.5) circle (1pt);
		\draw[fill=black] (4.5,1.5) circle (1pt);
		\draw[fill=black] (4.5,2.5) circle (1pt);
		\draw[fill=black] (5.5,2.5) circle (1pt);
		\draw[fill=black] (5.5,3.5) circle (1pt);
		\draw[fill=black] (6.5,3.5) circle (1pt);
		\draw[fill=black] (6.5,4.5) circle (1pt);
		\draw[fill=black] (7.5,3.5) circle (1pt);
		\draw[fill=black] (7.5,4.5) circle (1pt);
		\draw[fill=black] (8.5,4.5) circle (1pt);
		
		\draw[ultra thick, Purple, rounded corners=1] (0,0) -- (1,0) -- (1,1) -- (3,1) -- (3,2) -- (5,2) -- (5,3) -- (6,3) -- (6,4) -- (8,4) -- (8,5) -- (8.8,5);
		
		\node[Purple] at (8.45,5.3) {$\hat\gamma_u$};
		\node[orange] at (9.3,4.825) {\small$\R^+u$};
		\end{tikzpicture}
	\end{center}
	Whenever $u=\atob$ with $a,b$ odd integers, the ray passes through the center of infinitely many squares. We choose that $\hat\gamma_{u}$ alternates above and below these squares.

    \medbreak
    
	Note that $\hat \gamma_u$ is geodesic in $\Z^2$, so $\gamma_u$ is an infinite geodesic in $\Cc$, hence defines a Busemann point in $\partial(\Cc,d_S)$. We will see these points are distinct in $\partial^r(\Cc,d_S)$.

\subsection{Computations of lengths}



We first prove a lower bound on the length of certain elements.
\begin{lemma} \label{lem:lower}
	If $g\in\Cc$ such that $\hat g=\hat\gamma_{u,n}$ and $\norm{g}_S=n+\Delta$, then
	\[ \la B(g) ; u^\perp \ra \le \la B(\gamma_{u,n}); u^\perp\ra  + O(\Delta^3) \]
	where $u^\perp$ is the image of $u$ under a rotation of $+90^\circ$. As a corollary, we have that
	\[ \norm{h\gamma_{u,n}}_S = n + \Omega\Big(\sqrt[3]{\langle B(h);u^\perp\rangle}\ \Big) \]
	for $h\in [\Cc,\Cc]$ with $\langle B(h);u^\perp\rangle>0$. 
\end{lemma}
\begin{remark}
    If $\langle B(h);u^\perp\rangle\le 0$, we only get the (trivial) inequality $\norm{h\gamma_{u,n}}\ge n$.
\end{remark}
The proof follows the same scheme as \cite[Proposition 4.2]{Intermediate}.

\begin{proof} We first explain how the first inequality implies the corollary. Take $g=h\gamma_{u,n}$. As $h\in[\Cc,\Cc]$, we indeed have $\hat g=\hat \gamma_{u,n}$ and
\[ B(h\gamma_{u,n}) = B(h) + B(\gamma_{u,n}) \]
so, if we define $\Delta=\norm{g}_S-n$, we have that
\[ \langle B(h);u^\perp \rangle \le O(\Delta^3) \iff \Delta=\Omega\Big(\sqrt[3]{\langle B(h);u^\perp\rangle}\Big).\]

\bigbreak

\noindent $\blacktriangleright$ Now, the main inequality! As $\hat g=\hat\gamma_{u,n}$, we have $B(g)-B(\gamma_{u,n})=B(gx)-B(\gamma_{u,n}x)$ for all $x\in\Cc$, so we are going to estimate $\langle B(g\delta);u^\perp\rangle$ instead of $\langle B(g);u^\perp\rangle$, where $\delta$ is the interval of $\gamma_u$ between $\hat\gamma_{u,n}$ and the next intersection with the ray $\R^+u$.
\begin{center}
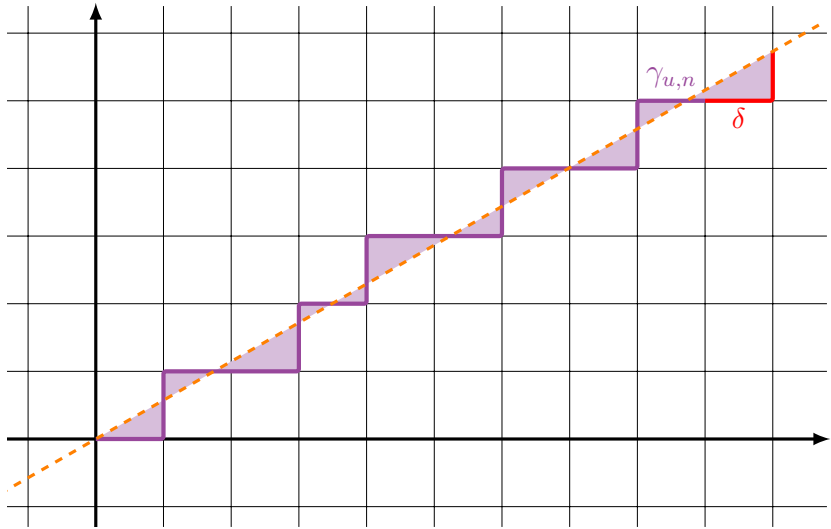

	\begin{tikzpicture}[scale=.9]
	\clip (-1.3,-1.3) rectangle (10.85,6.45);
	
	\draw[very thin] (-2.3,-1.3) grid (10.8,6.4);
	
	\draw[-latex, very thick] (-2.35,0) -- (10.85,0);
	\draw[-latex, very thick] (0,-1.35) -- (0,6.45);		
	
	\fill[Purple!30, opacity=.5, rounded corners=1] (0,0) -- (1,0) -- (1,1) -- (3,1) -- (3,2) -- (4,2) -- (4,3) -- (6,3) -- (6,4) -- (8,4) -- (8,5) -- (10,5) -- (10,5.73);
	\draw[ultra thick, Purple, rounded corners=1] (0,0) -- (1,0) -- (1,1) -- (3,1) -- (3,2) -- (4,2) -- (4,3) -- (6,3) -- (6,4) -- (8,4) -- (8,5) -- (9,5);
	\draw[ultra thick, red, rounded corners=1] (9,5) -- (10,5) -- (10,5.73);
	
	\node[Purple] at (8.5,5.33) {$\gamma_{u,n}$};
	\node[red] at (9.5,4.75) {\small$\delta$};
	
	\draw[orange, dashed, very thick] (-2.225,-1.275) -- (10.68,6.12); 
	\end{tikzpicture}
    \captionof{figure}{The path $\delta_{u,n}\delta$}
\end{center}

Fix a geodesic path representing $g$. We decompose the concatenation $g\delta$ as follows: we first remove all loops, and draw them in blue. We also cut the parts until the last crossing of $\R^-u$, and from the first crossing of \say{$\R^+u$ after $\widehat{g\delta}$}. The remaining is the \say{main part}, drawn in green. The green simple curve and the blue multi-curve both define a winding number distribution (see Figure \ref{fig:decomposition}), hence both of them have a well-defined (non-normalized) barycenter. Moreover, we have
\[ B(g\delta)=B(\text{green area})+B(\text{blue area}), \]
hence we can estimate both contributions to $\langle B(g\delta);u^\perp\rangle$ separately.

\begin{center}
	\begin{tikzpicture}[scale=.9]
\begin{scope}[shift={(0,0)}]
    \clip (-1.3,-1.3) rectangle (10.85,6.45);
	\draw[thin] (-2.3,-1.3) grid (10.8,6.4);
	
	\draw[-latex, very thick] (-2.35,0) -- (10.85,0);
	\draw[-latex, very thick] (0,-1.35) -- (0,6.45);		
	
	\fill[Green!30, opacity=.5, rounded corners=1] (-1,-.57) -- (-1,2) -- (5,2) -- (5,4) -- (9,4) -- (9,5) -- (10,5) -- (10,5.73);
	\fill[Turquoise!30, opacity=.5, rounded corners=1] (0,0) -- (0,1) -- (1,1) -- (1,-1) -- (-1,-1) -- (-1,-.57)
	(5,4) -- (5,6) -- (3,6) -- (3,4) -- (5,4);
	\draw[ultra thick, Turquoise, rounded corners=1] (0,0) -- (0,1) -- (1,1) -- (1,-1) -- (-1,-1) -- (-1,-.57)
	(5,4) -- (5,6) -- (3,6) -- (3,4) -- (5,4);
	\draw[line width=2pt, Turquoise, -latex] (3,5.15) -- (3,4.75);
	\draw[ultra thick, Green, rounded corners=1] (-1,-.57) -- (-1,2) -- (5,2) -- (5,4) -- (9,4) -- (9,5);
	
	\draw[ultra thick, red, rounded corners=1] (9,5) -- (10,5) -- (10,5.73);
	
	\node[Green] at (8.5,4.45) {\small$+1$};
	\node[Green] at (5.45,3.6) {\small$-1$};
	\node[Green] at (4.55,2.3) {\small$+1$};
	\node[Green] at (.5,1.5) {\small$-1$};
	\node[Turquoise] at (.5,-.5) {\small$-1$};
	
	\draw[orange, dashed, very thick] (-2.225,-1.275) -- (10.68,6.12); 
	\node[Green] at (8.5,3.7) {$g$};
	\node[red] at (9.5,4.75) {\small$\delta$};
\end{scope}

\begin{scope}[shift={(-2,-6)}, scale=.6]
    \draw[-latex, very thick] (-2,0) -- (10.85,0);
	\draw[-latex, very thick] (0,-1.35) -- (0,6.45);		
	
	\fill[Green!30, opacity=.5, rounded corners=1] (-1,-.57) -- (-1,2) -- (5,2) -- (5,4) -- (9,4) -- (9,5) -- (10,5) -- (10,5.73); 
    \draw[ultra thick, Green, rounded corners=1, -latex] (-1,-.57) -- (-1,2) -- (5,2) -- (5,4) -- (9,4) -- (9,5) -- (10,5) -- (10,5.73);
    \draw[thick, dashed, orange] (-1,-.57) -- (10,5.73);
    \node[Green] at (8.5,4.45) {\scriptsize$+1$};
	\node[Green] at (5.45,3.6) {\scriptsize$-1$};
	\node[Green] at (4.55,2.3) {\scriptsize$+1$};
	\node[Green] at (.7,1.2) {\footnotesize$-1$};
\end{scope}

\begin{scope}[shift={(7,-6)}, scale=.6]
    \draw[-latex, very thick] (-2,0) -- (10.85,0);
	\draw[-latex, very thick] (0,-1.35) -- (0,6.45);
 
	\draw[dotted, Green, rounded corners=1] (-1,-.57) -- (-1,2) -- (5,2) -- (5,4) -- (9,4) -- (9,5) -- (10,5) -- (10,5.73);
	\fill[Turquoise!30, opacity=.5, rounded corners=1] (0,0) -- (0,1) -- (1,1) -- (1,-1) -- (-1,-1) -- (-1,-.57)
	(5,4) -- (5,6) -- (3,6) -- (3,4) -- (5,4);
	\draw[ultra thick, Turquoise, rounded corners=1] (0,0) -- (0,1) -- (1,1) -- (1,-1) -- (-1,-1) -- (-1,-.57)
	(5,4) -- (5,6) -- (3,6) -- (3,4) -- (5,4);
	\draw[line width=2pt, Turquoise, -latex] (3,5.15) -- (3,4.75);
    \draw[line width=2pt, Turquoise, -latex] (1,0) -- (1,-.4);
    \draw[thick, dashed, Turquoise] (0,0) -- (-1,-.57);

    \node[Turquoise] at (.5,-.5) {\scriptsize$-1$};
    \node[Turquoise] at (4,5) {\footnotesize$+1$};
\end{scope}
\end{tikzpicture}
    
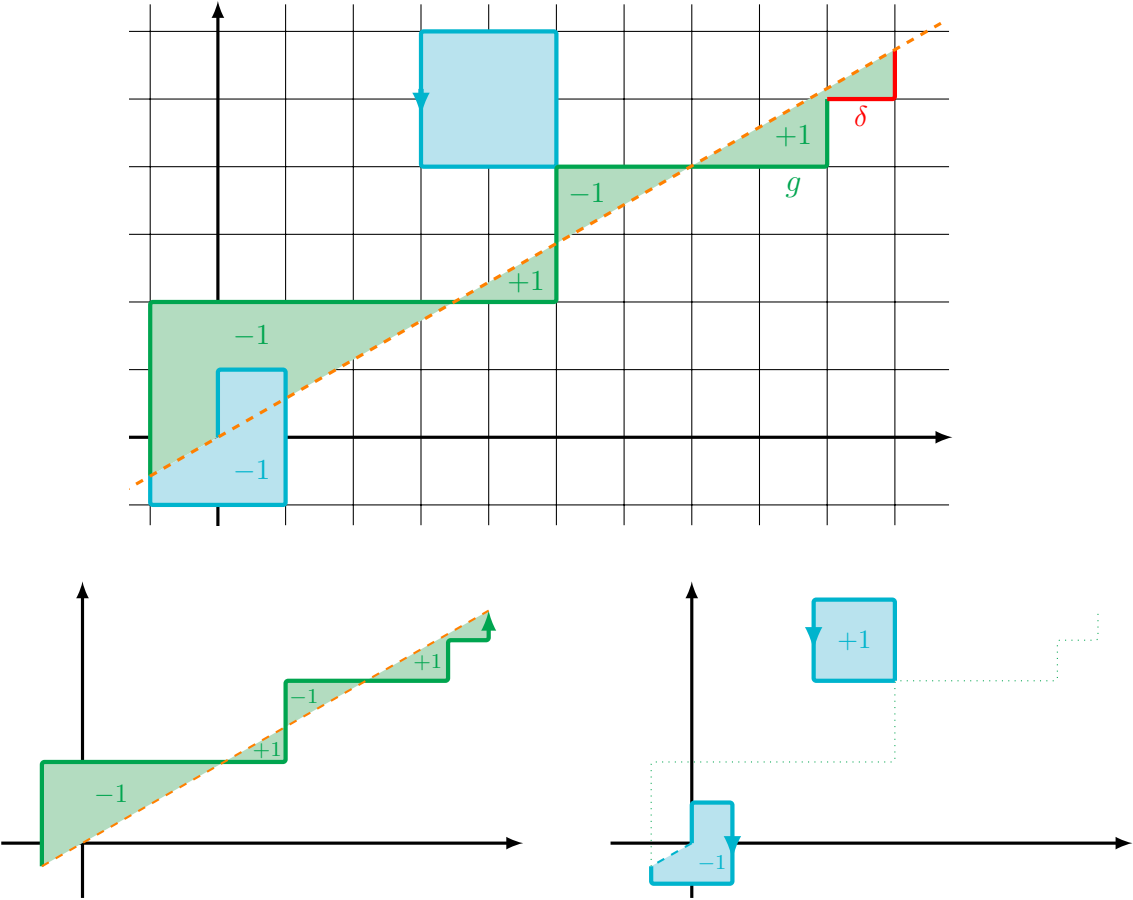
\captionof{figure}{Decomposition of geodesic representing $g$} \label{fig:decomposition}
\end{center}

\medskip

\noindent$\blacktriangleright$ Let's first suppose the entire path $g\delta$ stays within a $2\Delta$-neighborhood of $\R u$.
\medbreak
First observe that, by construction, both $\gamma_{u,n}\delta$ and the green path are simple curves. Their winding numbers (after closing using $\R u$) are $\pm 1$ or $0$. More specifically, $+1$'s only appear below $\R u$, and $-1$'s only appear above $\R u$, so the contribution of each region between $g\delta$ (or $\gamma_{u,n}\delta$) and $\R u$ to $\langle B;u^\perp\rangle$ is negative. Moreover, for each square cut in two by $\R u$, one of the two halves has to have non-zero winding number. Now we realize that $\gamma_{u,n}\delta$ is defined so that $\langle B(\gamma_{u,n}\delta);u^\perp\rangle$ is closest to $0$ among simple curves reaching the same endpoint, we have
\[ \la B(\text{green area});u^\perp\ra \le \la B(\gamma_{u,n}\delta);u^\perp\ra<0. \]
On the other hand, for the blue area, the length of the blue loops is bounded by $2\Delta$. It follows that the total blue area is bounded by $I(2\Delta)^2$ for some isoperimetric constant $I>0$ (for the grid, $I=\frac1{16}$). As all this area lies within $2\Delta$ from $\R u$, we get that
\[\la B(\text{blue area}); u^\perp \ra \le 8I \Delta^3 \]
and finally $\la B(g\delta);u^\perp\ra \le \la B(\gamma_{u,n}\delta);u^\perp\ra + 8I \Delta^3$.

\bigbreak

\noindent$\blacktriangleright$ Let's suppose the furthest $g\delta$ gets from $\R u$ is exactly $L\Delta$, for some $L\ge 2$.

\medbreak

As before, the blue contribution is bounded by $\la B(\text{blue area}); u^\perp \ra \le 4L I \Delta^3$. On the other hand, the curve $g\delta$ goes through a point $p$ at distance $L\Delta$ from $\R u$. This forces the green area to contain a large triangle, disjoint from the forced half-cut squares.
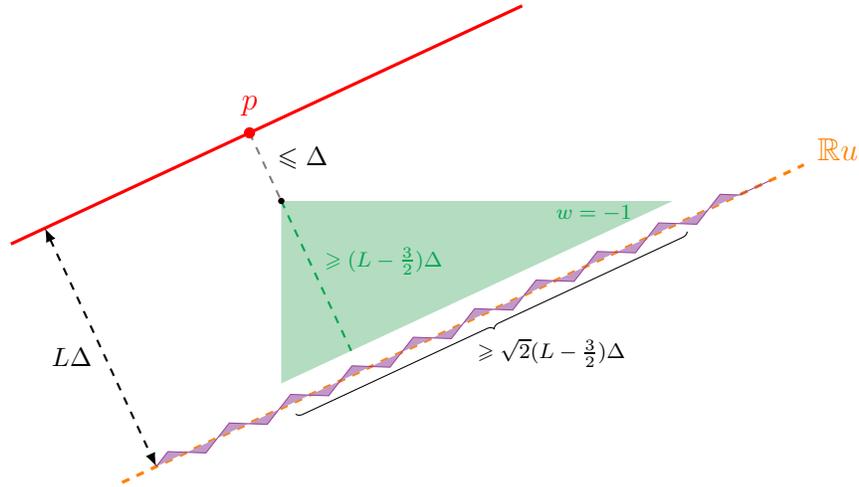
\begin{center}
	\begin{tikzpicture}[scale=1, rotate=25]
	\fill[Green!30, opacity=.5] (3.5,2.5) -- ({3.5-2.2*tan(25)},.3) -- ({3.5+2.2*tan(65)},.3);
	
	\draw[orange, dashed, very thick] (0,0) -- (10,0);
	\draw[red, very thick] (0,3.5) -- (7.5,3.5);
	
	\draw[thick, dashed, gray] (3.5,3.5) -- (3.5,2.5);
	\draw[thick, dashed, Green] (3.5,2.5) -- (3.5,.3);
	\draw[thick, dashed, latex-latex] (0.5,3.5) -- (0.5,0);
	\draw[scale=2, Purple, decorate, decoration={zigzag, segment length=9mm, amplitude=1.15mm}, fill=Purple!50, opacity=.5] (.25,0) -- (4.75,0);
	\node[circle, inner sep=1.5pt, fill=red, label=above:$\color{red}p$] at (3.5,3.5) {};
	\draw[fill=black] (3.5,2.5) circle (1pt);
	\node at (0.1,1.8) {\footnotesize$L\Delta$};
	\node at (4,2.9) {\footnotesize$\le\Delta$};
	\node[Green] at (4.4,1.2) {\scriptsize$\ge (L-\frac32)\Delta$};
	\node[orange] at (10.5,0) {$\R u$};
	\draw[decorate, decoration=brace]  ({3.5+2.2*tan(65)},-.15) -- ({3.5-2.2*tan(25)},-.15)
	node [black,midway, xshift=8mm, yshift=-4mm] {\scriptsize$\ge \sqrt2(L-\frac32)\Delta$};
	
	\node[Green] at (7.2,.6) {\scriptsize$w=-1$};
	\end{tikzpicture}
\end{center}
It follows that
\begin{align*}
\la B(\text{green area});u^\perp)\ra
& \le \la B(\gamma_{u,n}\delta);u^\perp\ra + \la B(\text{green triangle});u^\perp\ra \\
& \le \la B(\gamma_{u,n}\delta);u^\perp\ra -\frac{\sqrt 2}6\left(L-\frac32\right)^3\Delta^3 
\end{align*}
hence $\langle B(g\delta);u^\perp\rangle \le \langle  B(\gamma_{u,n}\delta);u^\perp\rangle + \Big(4L I -\frac{\sqrt 2}6\left(L-\frac32\right)^3\Big)\Delta^3$.

\bigbreak

\noindent$\blacktriangleright$ Finally, we observe that $L\mapsto 4LI-\frac{\sqrt2}6\big(L-\frac32\big)^3$ is eventually decreasing, so
\[ \la B(g\delta);u^\perp\ra \le \la B(\gamma_{u,n}\delta);u^\perp\ra + M \Delta^3 \]
where $M=\max\big\{8I,\; 4LI-\frac{\sqrt2}6\big(L-\frac32\big)^3 \mid L\ge 2\big\}$.
\end{proof}

\bigbreak

Next we give an upper bound on the length of specific elements.

\begin{lemma}\label{lem:upper} For $h\in[\Cc,\Cc]$ and $n$ large enough {\normalfont(}depending on $h${\normalfont)}, we have 
	\[ \norm{h\gamma_{u,n}}_S = \begin{cases}
	n + O\big(\sqrt[3]{\la B(h);u^\perp\ra + \abs{A(h)}}\big) + O(1) & \text{if }\la B(h);u^\perp\ra \ge 0, \\
	n + O\big(\sqrt[3]{\abs{A(h)}}\big)+O(1) & \text{if }\la B(h);u^\perp\ra\le 0.
	\end{cases} \]
	If moreover $u$ is \textbf{not} of the form $\atob$ for integers $a,b$ such that $a$ xor $b$ is even, then we can improve the estimates and eliminate the dependence on $\abs{A(h)}$:
	\[ \norm{h\gamma_{u,n}}_S = \begin{cases} n+O\big(\sqrt[3]{\la B(h);u^\perp\ra}\big)+O(1) & \text{if }\la B(h);u^\perp\ra \ge 0, \\ n+O(1) & \text{if }\la B(h);u^\perp\ra\le 0. \end{cases} \]
\end{lemma}

\begin{proof} We are looking for a short path representing $h\gamma_{u,n}$. We modify the path $\gamma_{u,n}$ (with $n$ large) to change the values of $A$ and $B$ from $A(\gamma_{u,n})$ and $B(\gamma_{u,n})$ to
\[ A(h\gamma_{u,n}) = A(h) + A(\gamma_{u,n}) \quad\text{and}\quad B(h\gamma_{u,n}) = B(h) + B(\gamma_{u,n}). \]
We fine tune the parameters $A$, $\la B;u\ra$ and $\langle B;u^\perp\rangle$ in three steps at the appropriate cost under the different hypothesis. The different operations are local, hence the endpoint $\hat\gamma_{u,n}$ remains unchanged all along. We control how the length $\ell$ of the considered path evolves after each step. Initially, $\ell=\norm{\gamma_{u,n}}_S=n$. 

\bigbreak

\noindent$\blacktriangleright$ First, we modify $A$ by exactly $p$ ($p\in\Z$) at a cost $O(1)$. Let's first suppose $u\ne (\pm1,0),(0,\pm1)$. Suppose $p>0$. We find $p$ squares which intersect the ray $\R u$, around which $\hat\gamma_u$ goes clockwise.
		
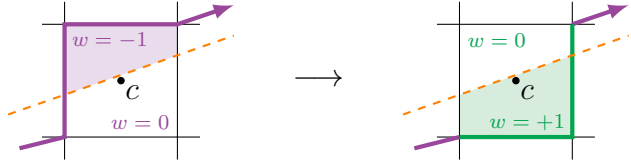
\begin{center}
			\begin{tikzpicture}[scale=1.5]
			
			\fill[Purple!15] (0,.375) -- (0,1) -- (1,1) -- (1,.725);
			\draw (-.2,-.2) grid (1.2,1.2);
			\draw[ultra thick, Purple, -latex] (-.4,-.1) -- (0,0) -- (0,1) -- (1,1) -- (1.5,1.17);
			
			\draw[thick, dashed, orange] (-.5,.2) -- (1.5,.9);
			\draw[fill=black] (.5,.5) circle (.8pt);
			\node at (.6,.4) {$c$};
			\node[Purple] at (.4,.86) {\scriptsize$w=-1$};
			\node[Purple] at (.67,.13) {\scriptsize$w=0$};
			
			\node at (2.25,.5) {$\longto$};
			
			\begin{scope}[shift={(3.5,0)}]
			\fill[Green!15] (0,.375) -- (0,0) -- (1,0) -- (1,.725);
			\draw (-.2,-.2) grid (1.2,1.2);
			\draw[ultra thick, Purple, -latex] (-.4,-.1) -- (0,0) (1,1) -- (1.5,1.17);
			\draw[ultra thick, Green] (0,0) -- (1,0) -- (1,1);
			
			\draw[thick, dashed, orange] (-.5,.2) -- (1.5,.9);
			\draw[fill=black] (.5,.5) circle (.8pt);
			\node at (.6,.4) {$c$};
			\node[Green] at (.33,.86) {\scriptsize$w=0$};
			\node[Green] at (.6,.13) {\scriptsize$w=+1$};
			\end{scope}
			\end{tikzpicture}
		\end{center}
		Each time we flip a unit square (from clockwise to counter-clockwise), we add $1$ to the winding numbers on the full square, so the parameters change by
		\begin{align*}
		A & \to A+1, \\
		B & \to B+c
		\end{align*}
		The distance from $c$ to $\R u$ being at most $\sqrt 2/2$, the component $\langle B;u^\perp\rangle$ changes by at most $\sqrt2/2$. On the other hand $\la B;u\ra$ changes essentially by $\norm{c}_E$. After flipping all $p$ squares, the total effect of this operation is
		\begin{align*}
		A & \to A+p, \\
		\la B;u\ra & \to \la B;u\ra + O(p^2) \\
		\langle B;u^\perp\rangle & \to \langle B;u^\perp\rangle+O(p), \\
		\ell & \to \ell.
		\end{align*}
        (If we take the first $p$ squares along the line around which $\hat\gamma_u$ goes clockwise, that $\norm{c}_E=O(p)$ for each square, which explains why $\la B;u\ra\to \la B;u\ra + O(p^2)$. As we will see later, the size of this change actually doesn't matter.)
		If $u=(\pm1,0),(0,\pm1)$, this doesn't quite work, we instead do the following change:
		
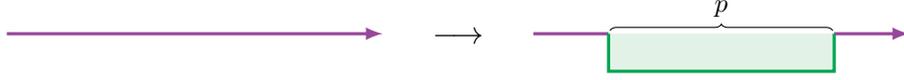
\begin{center}
			\begin{tikzpicture}
				\draw[very thick, Purple, -latex] (0,0) -- (5,0);
				
				\node at (6,-.05) {$\longto$};
				
				\begin{scope}[shift={(7,0)}]
					\draw[very thick, Purple, -latex] (0,0) -- (1,0) (4,0) -- (5,0);
					\draw[very thick, Green, fill=Green!10] (1,0) -- (1,-.5) -- (4,-.5) -- (4,0);
					\draw[decorate, decoration=brace] (1.02,.04) -- (3.98,.04)
					node [black,midway, yshift=3mm] {\footnotesize$p$};
				\end{scope}
			\end{tikzpicture}
		\end{center}
		with the same effect except $\ell\to\ell+2$.
  
	\bigbreak
  
	\noindent\textbf{Under the hypothesis} that $u$ is not of the form $\atob$ with $a,b$ integers such that $a$ xor $b$ is even, we can do better. We modify $A$ by exactly $p$, at cost $0$, and only modifying $\langle B;u^\perp\rangle$ by $O(1)$ (and not $O(p)$ as before).

    \medbreak
    
	Thanks to the hypothesis, we can find unit squares with center $c$ either on the ray $\R^+u$ (if $u=\atob$ with $a,b$ odd integers), or at distance at most $\frac1{\abs p}$ from the ray (for irrational slopes). Then, as before, we can select $\abs p$ of them to flip from clockwise to counter-clockwise (if $p>0$) or the other way around (if $p<0$). After flipping all $p$ squares, the total effect of this operation is
		\begin{align*}
		A & \to A+p, \\
		\la B;u\ra & \to \la B;u\ra + O_u(p^3) \\
		\langle B;u^\perp\rangle & \to \langle B;u^\perp\rangle+O(1), \\
		\ell & \to \ell.
		\end{align*}
		(The variation in $\la B;u\ra$ may be huge, but this is nothing to worry about.\footnote{The subscript $u$ indicates the implied constants in $O_u(p^3)$ might depend on $u$.})
		
		\bigbreak
		
\noindent$\blacktriangleright$ Next we change $\langle B;u^\perp\rangle$ by approximately $q$ at the appropriate cost.
\begin{itemize}[leftmargin=8mm, label=$\star$]
			\item If $q>0$, we can find $b\in\Z^2$ such that $\langle b;u^\perp\rangle=q+O(1)$ and $\la b;u\ra =O(1)$ (hence $\norm{b}_E=q+O(1)$), and $z\in[\Cc,[\Cc,\Cc]]$ such that $B(z)=b$.
			
			Because $[\Cc,[\Cc,\Cc]]$ is cubically distorted inside $\Cc$, we have $\norm{z}_S=O(\sqrt[3]{q})$. (See for instance \cite[Corollary 14.16]{GGT}.) So we can find a loop of this length evaluating to $z$, and then glue it at any point along $\gamma_u$. The total effect is
			\begin{align*}
			A & \to A \\
			B & \to B + b \\
			\la B;u\ra & \to \la B;u\ra + O(1) \\
			\langle B;u^\perp\rangle & \to \langle B;u^\perp\rangle + q + O(1) \\
			\ell & \to \ell+O(\sqrt[3]q)
			\end{align*}\vspace*{3mm}

			\item \begin{minipage}[t]{.85\linewidth}
				If $q<0$. There exists $v\in\{(\pm1,0)(0,\pm1)\}$ s.t.\ $\langle v;u\rangle\ge\sqrt2/2$, without lost of generality let's assume that $v=(1,0)$. In  particular we also have $\langle (0,1);u^\perp\rangle\ge \sqrt 2/2$. \vspace{2mm}
			\end{minipage} \hspace*{3mm}
			\begin{minipage}[t]{.10\linewidth}
				\vspace{-8mm}
				\begin{tikzpicture}[scale=.7]
				\draw (0,0) circle (1cm);
				\draw[orange, ultra thick, domain=-45:45] plot ({cos(\x)}, {sin(\x)});
				\draw[LimeGreen, ultra thick, domain=45:135] plot ({cos(\x)}, {sin(\x)});
				\node[circle, thick, draw=black, fill=orange, inner sep=1.2pt, label=right:$u$] at ({cos(25)}, {sin(25)}) {};
				\node[circle, thick, draw=black, fill=LimeGreen, inner sep=1.2pt] at ({cos(115)}, {sin(115)}) {};
				\node at (-.6,1.2) {$u^\perp$};
				\end{tikzpicture}
			\end{minipage}
		
			For $n$ large, we can find two disjoint sub-strings $\alpha,\beta$ in the untouched part of the path $\hat\gamma_u$, with $\abs{q}$ occurrences of the letter $x$ each (because $\la (1,0);u\ra>0$). We assume that both $\alpha,\beta$ start with an $x$.
			
			For $0\le k\le \abs q$, let $\alpha_k$ be the prefix of $\alpha$ which stop right after the $k^\text{th}$ occurrence of $x$. By convention $\alpha_0=\emptyset$. We define $\beta_k$ similarly. The operation is the following: we replace the string $\alpha_k$ by $y^2\alpha_ky^{-2}$, and the string $\beta_k$ by $y^{-2}\beta_ky^2$. We prove that, for an appropriate $k$, the value $\langle B;u^\perp\rangle$ changes by $q+O(1)$.
   \begin{adjustwidth}{-15mm}{-15mm}
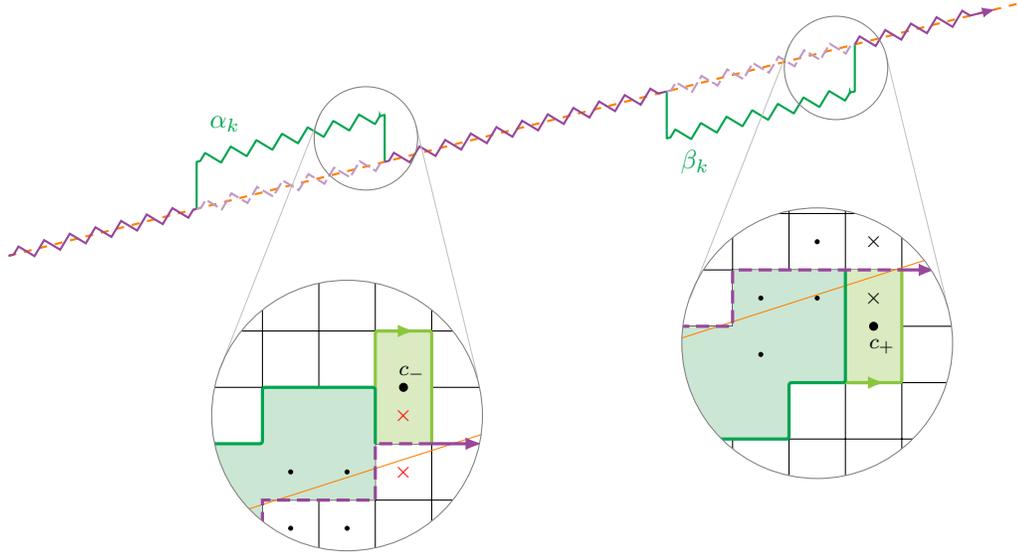
	\begin{center}
		\begin{tikzpicture}[scale=1.25]
	
	\draw[dashed, orange, thick] (0,0) -- (10.75,2.6875);
	\draw[thick, Purple, decorate, decoration={zigzag, pre length=.7mm}] (0,0)  -- (2,.5);

	\draw[thick, Green, rounded corners=.2] (2,.5) -- (2,1)  
	decorate [decoration={zigzag, pre length=.5mm, post length=.5mm}]
	{to  (4,1.5)} -- (4,1);;
	\draw[thick, Purple!50, dash pattern=on 4.5pt off 1.5pt, decorate, decoration={zigzag, pre length=.6mm}]
	(2,.5)  -- (4,1);
	\draw[thick, Purple, decorate, decoration={zigzag, pre length=.6mm}] (4,1)  -- (7,1.75);
	\draw[thick, Green, rounded corners=.2] (7,1.75) -- (7,1.25)
	decorate [decoration={zigzag, pre length=.5mm, post length=.5mm}] 
	{to (9,1.75)}  -- (9,2.25);
	\draw[thick, Purple!50, dash pattern=on 4.5pt off 1.5pt, decorate, decoration={zigzag, pre length=.6mm}] (7,1.75)  -- (9,2.25);
	\draw[thick, Purple, decorate, decoration={zigzag, post length=2mm}, -latex] (9,2.25)  -- (10.5,2.625);
	
	\node[Green] at (2.3,1.4) {\footnotesize$\alpha_k$};
	\node[Green] at (7.3,1) {\footnotesize$\beta_k$};
	
	\draw[very thin, lightgray] (2.25,-1.2) -- (3.31,1.5);
	\draw[very thin, lightgray] (4.99,-1.3) -- (4.335,1.38);
	
	\draw[very thin, lightgray] (7.25,-.45) -- (8.31,2.25);
	\draw[very thin, lightgray] (9.99,-.55) -- (9.335,2.13);
	
	\draw[gray, thin] (3.8,1.25) circle (.55cm);
	\draw[gray, thin] (8.8,2) circle (.55cm);
	
	\begin{scope}[shift={(2.7,-2.6)}, scale=.6]
		\begin{scope}
			\clip (1.5,1.5) circle (2.4cm);
			\draw (-1,-1) grid (4,4);
			\fill[Green!20, opacity=.5] (-1,-1) -- (0,-1) -- (0,0) -- (2,0) -- (2,2) -- (0,2) -- (0,1) -- (-1,1) -- (-1,-1);
			\fill[LimeGreen!30, opacity=.5] (2,1) -- (2,3) -- (3,3) -- (3,1) -- (2,1);
			
			\draw[orange] (-1,-.4) -- (4,1.2);
			\draw[very thick, LimeGreen, rounded corners=1] (2,1.96) -- (2,3) -- (3,3) -- (3,1);
			\draw[very thick, LimeGreen, -latex] (2.1,3) -- (2.7,3);
			\draw[very thick, Green, rounded corners=1] (-1,1) -- (0,1) -- (0,2) -- (2,2) -- (2,1);
			\draw[very thick, dash pattern=on 5pt off 2.5pt, Purple] (0,-1.5) -- (0,0) -- (2,0) -- (2,1) -- (3,1);
			\draw[very thick, Purple, -latex] (3,1) -- (3.9,1);
			
			\draw[fill=black] (.5,.5) circle (1pt);
			\draw[fill=black] (.5,-.5) circle (1pt);
			\draw[fill=black] (1.5,.5) circle (1pt);
			\draw[fill=black] (1.5,-.5) circle (1pt);
			\node[red] at (2.5,0.5) {\scriptsize$\times$};
			\node[red] at (2.5,1.5) {\scriptsize$\times$};
			
			\draw[fill=black] (2.5,2) circle (2pt);
			\node at (2.65,2.25) {\scriptsize$c_-$};
		\end{scope}
		\draw[gray] (1.5,1.5) circle (2.4cm);
	\end{scope}
	
	\begin{scope}[shift={(7.7,-1.95)}, scale=.6]
		\begin{scope}
			\clip (1.5,1.7) circle (2.4cm);
			\draw (-1,-1) grid (4,4.5);
			\fill[Green!20, opacity=.5] (-1,0)--(1,0)--(1,1) -- (2,1) -- (2,3) -- (0,3) -- (0,2) -- (-1,2) -- (-1,0);
			\fill[LimeGreen!30, opacity=.5] (2,1) -- (2,3) -- (3,3) -- (3,1) -- (2,1);
	
			\draw[orange] (-1.5,1.6) -- (3.5,3.2);
			\draw[very thick, LimeGreen, rounded corners=1] (1.96,1) -- (3,1) -- (3,3);
			\draw[very thick, Green, rounded corners=1] (-1,0) -- (1,0) -- (1,1) -- (2,1) -- (2,3);
			\draw[very thick, LimeGreen, -latex] (2.1,1) -- (2.7,1);
			\draw[very thick, dash pattern=on 5pt off 2.5pt, Purple] (-1,2) -- (0,2) -- (0,3) -- (3,3);
			\draw[very thick, Purple, -latex] (3,3) -- (3.6,3);
	
			\draw[fill=black] (.5,1.5) circle (1pt);
			\draw[fill=black] (.5,2.5) circle (1pt);
			\draw[fill=black] (1.5,2.5) circle (1pt);
			\draw[fill=black] (1.5,3.5) circle (1pt);
			\node at (2.5,2.5) {\scriptsize$\times$};
			\node at (2.5,3.5) {\scriptsize$\times$};
	
			\draw[fill=black] (2.5,2) circle (2pt);
			\node at (2.65,1.65) {\scriptsize$c_+$};
		\end{scope}
		\draw[gray] (1.5,1.7) circle (2.4cm);
	\end{scope}

\end{tikzpicture}
	\end{center}
    \end{adjustwidth}
	What happens when we replace $k$ by $k+1$? The winding numbers inside two $1\times 2$ rectangles change by $-1$ and $+1$ respectively. So the effect is
	\begin{align*}
		A_k & \to A_{k+1} = A_k +2-2\\
		B_k & \to B_{k+1} = B_k + 2 (c_+-c_-)
	\end{align*}
	Recall that, by construction, the ray $\R^+u$ passes in between the two red crosses in the previous picture. In particular, 
		\[ \frac{\sqrt2}4 \le \la (0,0.5);u^\perp\ra \le \la c_-;u^\perp\ra \le \la(0,1.5);u^\perp\ra \le \frac32 \]
	Similarly, we get $-\frac32\le \langle c_+;u^\perp\rangle\le -\frac{\sqrt2}4$. Overall, this gives
		\[ -6 \le \la B_{k+1};u^\perp\ra - \la B_k;u^\perp \ra \le -\sqrt2 < -1\]
	It follows that we can find $k$ such that $\langle B_k;u^\perp\rangle-\langle B;u^\perp\rangle = \langle B_k;u^\perp\rangle-\langle B_0;u^\perp\rangle$ falls at distance at most $3$ from $q$. The total effect is
	\begin{align*}
		A & \to A \\
		\la B;u\ra & \to \la B;u\ra + O(q^2) \\
		\langle B;u^\perp\rangle & \to \langle B_k;u^\perp\rangle = \langle B;u^\perp\rangle + q + O(1) \\
		\ell & \to \ell+8
	\end{align*}
\end{itemize}	
		
\medbreak
		
\noindent$\blacktriangleright$ Finally we change $\la B;u\ra$ by approximately $r$ at a cost $O(1)$. We glue two loops $[x,y]^{-1}$ and $[x,y]$ on $\hat\gamma_u$, at a distance approximately $r$ from each other.
		
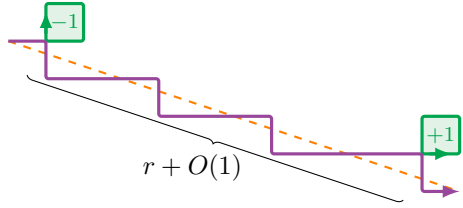
\begin{center}
			\begin{tikzpicture}[scale=.5]
				\draw[thick, dashed, orange] (0,0) -- (12,-4);

				\draw[very thick, Green, fill=Green!10, rounded corners=1] (1,0) -- (2,0) -- (2,1) -- (1,1) -- (1,0);
				\draw[very thick, Green, -latex] (1,.2) -- (1,.8);
				\draw[very thick, Green, fill=Green!10, rounded corners=1] (11,-3) -- (11,-2) -- (12,-2) -- (12,-3) -- (11,-3);
				\draw[very thick, Green, -latex] (11.2,-3) -- (11.8,-3);
				\node[Green] at (1.5,.5) {\scriptsize$-1$};
				\node[Green] at (11.5,-2.5) {\scriptsize$+1$};
				
				\draw[very thick, Purple, -latex, rounded corners=1] (0,0) -- (1,0) -- (1,-1) -- (4,-1) -- (4,-2) -- (7,-2) -- (7,-3) --( 11,-3) -- (11,-4) -- (12,-4);
				
				\draw[decorate, decoration=brace] (10.5,-4.23) -- (.5,-.9)
				node [black,midway,xshift=-3mm, yshift=-4mm] {\footnotesize$r+O(1)$};
			\end{tikzpicture}
		\end{center}
		Specifically, if $c_+$ and $c_-$ are the centers of the two square loops, we ensure that $\norm{(c_+ \hspace{.3mm}-\hspace{.5mm} c_-)-qu}_E=O(1)$. The total effect of this operation is then
		\begin{align*}
			A & \to A, \\
			B & \to B + (c_+ \hspace{.3mm}-\hspace{.5mm} c_-) \\
			\la B;u\ra & \to\la B;u\ra +r+O(1), \\
			\langle B;u^\perp\rangle & \to \langle B;u^\perp\rangle+O(1), \\
			\ell & \to \ell+8.
		\end{align*}
	
\noindent$\blacktriangleright$ Applying the three steps (in that order), we transform the curve $(\hat\gamma_{u,k})_{0\le k\le n}$ with parameters $A(\gamma_{u,n})$ and $B(\gamma_{u,n})$ into another curve with parameters
\[ A(h) + A(\gamma_{u,n})+O(1)\quad\text{and}\quad B(h)+B(\gamma_{u,n})+O(1). \]
For the second step, one need to take $q=\langle B(h);u^\perp\rangle+O(\abs{A(h)})$ in the general case, or $q=\langle B(h);u^\perp\rangle+O(1)$ under the extra hypothesis. This new curve has length $\ell=n+O(1)+O(\max\{\sqrt[3]q;1\})+O(1)$. Finally, we can fix the last $O(1)$ difference of the endpoint in $\Cc$ at an extra $O(1)$ cost. Overall, this gives the desired upper bound \[ \norm{h\gamma_{u,n}}_S\le \ell =n+O(\max\{\sqrt[3]q;1\}). \qedhere\]
\end{proof}

\bigbreak
\subsection{Back to horofunctions}

\begin{theorem} \label{prop:Bus_on_comms} There exists $C_1,C_2>0$ such that, for all $h\in[\Cc,\Cc]$, we have
	\begin{align*}
	C_1\sqrt[3]{\la -B(h);u^\perp\ra} & \le b_{\gamma_u}(h) \le C_2\sqrt[3]{\la -B(h);u^\perp\ra + \abs{A(h)}}+C_2 & \text{if }\la -B(h);u^\perp\ra \ge 0, \\
	0 & \le b_{\gamma_u}(h) \le C_2\sqrt[3]{\abs{A(h)}}+C_2 & \text{if }\la -B(h);u^\perp\ra\le 0.
	\end{align*}
	If moreover $u$ is \textbf{not} of the form $\atob$ for integers $a,b$ such that $a$ xor $b$ is even, then we can improve the estimates and remove the dependence on $A(h)\colon$
	\[b_{\gamma_u}(h) = \begin{cases} \Theta\big(\sqrt[3]{\la -B(h);u^\perp\ra}\big) & \text{if }\la -B(h);u^\perp\ra> 0, \\ O(1) & \text{if }\la -B(h);u^\perp\ra\le 0. \end{cases} \]
\end{theorem}
\begin{proof}
	This follows directly from Lemma \ref{lem:lower} and \ref{lem:upper}, given that
	\[ b_{\gamma_u}(h) = \lim_{n\to\infty}  d(\gamma_n,h)-d(\gamma_n,e)= \lim_{n\to\infty} \norm{h^{-1}\gamma_n}_S - n \]
	The only missing piece is $\norm{h^{-1}\gamma_{u,n}}_S\ge n$, which is true because $h^{-1}\gamma_{u,n}$ has length $n$ in the abelianization (same endpoint as $\gamma_{u,n}$).
\end{proof}
\begin{remark}
This should be compared with the discrete Heisenberg group $H_3(\Z)$, and more generally pairs $(G,S)$ with \say{property EH} of \cite{Bader}, for which $\varphi(h)=O(1)$ for all horofunction $\varphi\in\partial(G,S)$ and all $h\in [G,G]$.
\end{remark}
\bigbreak

Finally, we are ready to prove our main theorem
\begin{theorem}\label{thm:cartan}
All Busemann points $[b_{\gamma_u}]$ are distinct in the \emph{reduced} boundary $\partial^r(\Cc,d_S)$. Moreover, the stabilizer of $[b_{\gamma_u}]$ for the action $\Cc\acts \partial^r(\Cc,d_S)$ satisfies
\[ \Stab_\Cc([b_{\gamma_u}]) \le \{g\in\Cc:\hat g\in \R u\} \]
{\normalfont(}which is $[\Cc,\Cc]$ for $u$ irrational{\normalfont)}.
\end{theorem}
\begin{proof}
	Consider two \emph{distinct} directions $u,v\in\SS^1$. There exists $b\in \Z^2$ such that $\la -b;u^\perp\ra> 0$ and $\la -b;v^\perp\ra\le 0$. Consider $h\in[\Cc,[\Cc,\Cc]]$ such that $B(h)=b$, then
	\begin{align*}
		b_{\gamma_u}(h^n) & = \Theta(\sqrt[3]n) \\
		b_{\gamma_v}(h^n) & = O(1).
	\end{align*}
	In particular $b_{\gamma_u}-b_{\gamma_v}$ is not a bounded function, that is, $[b_{\gamma_u}]\ne[b_{\gamma_v}]$.
 
	\bigbreak
 
	\noindent\begin{minipage}{.7\linewidth}
	$\blacktriangleright$ Now let us have a look at $\Stab_\Cc([b_{\gamma_u}])$. Fix
		\begin{itemize}[leftmargin=6mm]
			\item $g\in \Cc$ such that $\hat g\notin \R u$ or equivalently $\langle \hat g;u^\perp\rangle\ne 0$
			\item $h\in [\Cc,\Cc]$ such that $A(h)>0$ and $B(h)=0$, for instance $h=[x,y][x^{-1},y^{-1}]$.
		\end{itemize}
	\end{minipage} \hspace*{10mm}
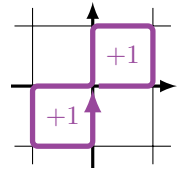
	\begin{minipage}{.25\linewidth}
		\centering
		\begin{tikzpicture}[scale=.8]
			\clip (-1.35,-1.35) rectangle (1.4,1.4);
			\draw (-1.3,-1.3) grid (1.3,1.3);
			
			\draw[-latex, very thick] (-1.4,0) -- (1.4,0);
			\draw[-latex, very thick] (0,-1.4) -- (0,1.4);
			
			\draw[-latex, Purple, line width=2pt, rounded corners=2] (0.1,0) -- (1,0) -- (1,1) -- (0,1) -- (0,0) -- (-1,0) -- (-1,-1) -- (0,-1) -- (0,-.05);
			
			\node[Purple] at (.5,.5) {\footnotesize$+1$};
			\node[Purple] at (-.5,-.5) {\footnotesize$+1$};
		\end{tikzpicture}
	\end{minipage}
 
	\bigbreak
 
	\noindent Let $C_1,C_2>0$ be the constants in Theorem \ref{prop:Bus_on_comms} such that
	\begin{align*}
	b_{\gamma_u}(*) & \ge C_1 \cdot \sqrt[3]{\langle -B(*);u^\perp\rangle} & \text{if }\langle -B(*);u^\perp\rangle\ge 0, \\
	b_{\gamma_u}(*) & \le C_2 \cdot \sqrt[3]{\abs{A(*)}} + C_2 & \text{if }\langle -B(*);u^\perp\rangle\le 0.
	\end{align*}
	Fix $m\in\Z$ such that $C_1\sqrt[3]{m\langle \hat g;u^\perp\rangle} > C_2$. (In particular, $m\langle \hat g;u^\perp\rangle>0$.) We can estimate the values of $\gamma_u$ and $g^m\cdot \gamma_u$ at $h^n$:
	\begin{align*}
	b_{\gamma_u}(h^n) & \le C_2 \cdot \sqrt[3]{nA(h)} + C_2 \\[4mm] 
	(g^m\cdot b_{\gamma_u})(h^n)\, & \hspace*{-1mm}:= b_{\gamma_u}(g^{-m}h^n)-b_{\gamma_u}(g^{-m}) \\[.5mm]
	& = b_{\gamma_u}\big(h^n[g^{-m},h^n]g^{-m}\big) + O_{g,u}(1) \\[.5mm]
	& = b_{\gamma_u}(h^n[g,h]^{-mn}) + O_{g,u}(1) \tag{1} \\
	& \ge C_1\sqrt[3]{m\langle \hat g;u^\perp\rangle} \cdot \sqrt[3]{nA(h)} + O_{g,u}(1) \tag{2}
	\end{align*}
	using that
 
    \medbreak
    
	\noindent (1) $b_{\gamma_u}$ is $1$-Lipschitz, we can extract the $g^{-m}$ from the right at a cost of $\norm{g^m}_S$.

    \medbreak
    
	\noindent (2) One can compute $B(h^n[g,h]^{-mn})=-mnB([g,h])=-mn\big(\hat g\cdot A(h)\big)$.

    \medbreak
    
	Letting $n\to+\infty$, we conclude that $b_{\gamma_u}$ and $g^m\cdot b_{\gamma_u}$ are not at bounded distance apart. This means $g^m$ doesn't fix $[b_{\gamma_u}]$, hence neither does $g$.
\end{proof}



\medskip

\begin{remark}
    The same computation shows that $\Ec\acts \partial^r(\Ec,d_S)$ is not trivial, for
    \[ \Ec = \la x,y \;\big|\; [x,[x,y]]=[x,[y,[x,y]]=[y,[y,[x,y]]=e \ra \]
    the Engel group (also called filiform group of step $3$), and $S=\{x^\pm,y^\pm\}$. More specifically, the orbit of $[b_{\gamma_u}]$ where $u=[1:0]$ (that is $[b_{x^\infty}]$) is infinite. However, it is not clear whether $\partial(\Ec,S)$ contains countably or uncountably many Busemann points.
\end{remark}

\medskip

Looking at larger groups, we get uncountably many Busemann points in many nilpotent groups, combining the previous result and Corollary \ref{lemma:quotient}.

\begin{corollary}
    Let $\G$ be a free nilpotent group of rank $r\ge 2$ and step $s\ge 3$ with the standard generating set $S$. Then $\partial^r(\G,d_S)$ contains uncountably many Buseman points.
\end{corollary}

\begin{proof}
    Since any free nilpotent group of rank $r\ge 2$ and step $s\ge 3$ surjects onto the Cartan group $\Cc$, the uncountably many Busemann points constructed in $\Cc$ lift to uncountably many Busemann points in $\G$, by Corollary \ref{lemma:quotient}.
\end{proof}

\printbibliography
\end{document}